\documentclass[12pt]{amsart}
\usepackage{amsthm,amssymb,amsmath,amscd}

\newtheorem{thm}{Theorem}[section]
\newtheorem{lm}[thm]{Lemma}
\newtheorem{cor}[thm]{Corollary}
\newtheorem{pro}[thm]{Proposition}
\newtheorem{pro-def}[thm]{Proposition and Definition}

\newtheorem{nota}[thm]{Notation}

\theoremstyle{definition}
\newtheorem{df}[thm]{Definition}

\input epsf

\def \int {\text{Int}}

\begin{document}

\title[Torsion of Quasi-Isomorphisms]
{Torsion of Quasi-Isomorphisms}
\author[J.-W. Chung and X.-S. Lin]{Jae-Wook Chung and Xiao-Song Lin}
\address{Department of Mathematics, University of California, Riverside, CA 92521}
\email{xl@math.ucr.edu, jwchung@alumni.ucr.edu}
\thanks{The second named author is supported in part by NSF}

\begin{abstract}{In this paper, we introduce the notion of Reidemeister
torsion for quasi-isomorphisms of based chain complexes over a
field. We call a chain map a quasi-isomorphism if its induced
homomorphism between homology is an isomorphism. Our notion of
torsion generalizes the torsion of acyclic based chain complexes,
and is a chain homotopy invariant on the collection of all
quasi-isomorphisms from a based chain complex to another. It
shares nice properties with torsion of acyclic based chain
complexes, like multiplicativity and duality. We will further
generalize our torsion to quasi-isomorphisms between free chain
complexes over a ring under some mild condition. We anticipate
that the study of torsion of quasi-isomorphisms will be fruitful
in many directions, and in particular, in the study of links in
3-manifolds.}
\end{abstract}

\maketitle

\section{Introduction}

The vector spaces used here are finite dimensional and rings are
commutative and have $1 \ne 0$.

It was observed by Milnor in his beautiful paper "Infinite cyclic
coverings" \cite{M} that the Alexander polynomial of a knot and
the Reidemeistor torsion of the infinite cyclic covering space of
the knot complement is directly related to each other, because of
the fact that the infinite cyclic covering of a knot complement is
acyclic when tensoring with the field of rational functions.

Turaev generalized this theorem of Milnor to the case of links
directly (see "Reidemeister torsion in knot theory" \cite{T2}).
But in the case of links, since the maximal abelian covering space
of a link complement can not be made acyclic in general, the
statement of Turaev's generalization is not as nice as that of the
theorem of Milnor.

One observation is that if we fix a link in $S^3$, say the trivial
link $L_0$, then there are infinitely many links $L$ in $S^3$,
which admit natural maps $S^3\setminus L\longrightarrow
S^3\setminus L_0$. These natural maps, when lifted to the
corresponding maximal abelian covering spaces, will induce
isomorphisms on homology after tensoring with the quotient field
of the polynomial ring over $\Bbb Z$. Can we extend the notion of
torsion to such a setting? If so, what will be the relationship
between such torsion defined in this setting and the Alexander
polynomial of $L$?

The goal of this paper is then to offer an approach for possibly
answering these questions. We introduce the notion of Reidemeister
torsion for quasi-isomorphisms of based chain complexes over a
field.  We call a chain map a quasi-isomorphism if its induced
homomorphism between homology is an isomorphism. Our notion of
torsion generalizes the torsion of acyclic based chain complexes,
and is a chain homotopy invariant on the collection of all
quasi-isomorphisms from a based chain complex to another. It
shares nice properties with torsion of acyclic based chain
complexes, like multiplicativity and duality. We will further
generalize our torsion to quasi-isomorphisms between free chain
complexes over a ring under some mild condition. Since the acyclic
condition is crucial whenever the notion of torsion is studied, we
hope that our point of view of torsion would be useful in other
directions.

The materials presented here contain only the most basic
definitions and proofs of basic properties of torsion of
quasi-isomorphisms. Although it was our original intention, we
have not yet worked out the details of the application of our
theory to the study of links in 3-manifolds. And the paper is
written in a somehow unsophisticated way, trying to cover as much
elementary details as possible. We hope that the reader will
tolerate us. The reader may also notice that our exposition
follows quite closely from that of the materials in Chapter 1 of
Turaev's book \cite{T}. Beside the introduction of the notion of
torsion of quasi-isomorphism, there are many technical details in
generalizing properties of torsion of acyclic chain complexes to
that of quasi-isomorphisms. We consider these as the main
contributions of this paper.

\section{Basic definitions and preliminaries}

In this section, we introduce basic terminologies and properties
used in this paper.

Let $V$ be a finite dimensional vector space over a field $F$.
Suppose that ${\rm dim}_FV=n$ and $b=(b_1,\dots,b_n)$ and
$b'=(b'_1,\dots,b'_n)$ are ordered bases for $V$. For convenience,
let us use row vectors. Then for each $i \in \{1,\dots,n\}$, there
is a unique $(a_{i1},\dots,a_{in}) \in F^n$ such that
$b_i=\sum_{j=1}^na_{ij}b'_j$, hence, we have the transition matrix
$(a_{ij})_{i,j=1,\dots,n}$ from $b$ to $b'$, denoted by $(b/b')$,
which is an $n \times n$ invertible matrix and write $[b/b']={\rm
det}(a_{ij})$.

\begin{pro} Define $\sim$ on the set $B$ of all ordered bases for
a finite dimensional vector space $V$ over a field $F$ by $b \sim
b'$ if and only if $[b/b']=1$ for all $b,b' \in B$. Then $\sim$ is
an equivalence relation on $B$. We call ordered bases $b$ and $b'$
equivalent if $b \sim b'$.
\end{pro}

\begin{proof} For each $b \in B$, $(b/b)$ is the identity matrix,
so $[b/b]=1$. That is, $b \sim b$. If $b \sim b'$ in $B$, then
$[b/b']=1$. Since $(b'/b)=(b/b')^{-1}$, $[b'/b]=1$. Hence, $b'
\sim b$. If $b \sim b'$ and $b' \sim b''$ in $B$, then $[b/b']=1$
and $[b'/b'']=1$. We claim that $(b/b'')=(b/b')(b'/b'')$. For each
$v \in V$, $v_{b'}=v_b(b/b')$ and $v_{b''}=v_{b'}(b'/b'')$. We
have $v_{b''}=v_b(b/b')(b'/b'')$. Hence, $(b/b'')=(b/b')(b'/b'')$,
so $[b/b'']=1$. Therefore, $b \sim b''$.
\end{proof}

\begin{lm} If $a$, $b$, $c$ are ordered bases for a finite
dimensional vector space $V$ over a field $F$, then $$a \sim b
\Leftrightarrow [a/c]=[b/c] \Leftrightarrow [c/a]=[c/b].$$
\end{lm}

\begin{proof} Notice that $[a/c]=[a/b][b/c]$. If $a \sim b$, then
$[a/b]=1$, so $[a/c]=[b/c]$. If $[a/c]=[b/c]$, then $[a/b]=1$, so
$a \sim b$. Also, since $[c/a]=[a/c]^{-1}$ and $[c/b]=[b/c]^{-1}$,
we have $[a/c]=[b/c] \Leftrightarrow [c/a]=[c/b]$.
\end{proof}

\begin{pro} If $F$ is a field and $A \in M_m(F)$, $B \in M_n(F)$,
$C \in M_{m \times n}(F)$, and $D \in M_{n \times m}(F)$, then
$${\rm det}\begin{pmatrix} A & C \\ 0 & B \end{pmatrix}={\rm
det}\begin{pmatrix} A & 0 \\ D & B \end{pmatrix}={\rm
det}\,A\,\,{\rm det}\,B.$$
\end{pro}

\begin{proof} By elementary row operations, we can change
$A$ and $B$ to upper triangular matrices $A'$ and $B'$,
respectively, so that ${\rm det}\,A=(-1)^r{\rm det}\,A'$ and ${\rm
det}\,B=(-1)^s{\rm det}\,B'$ for some nonnegative integers $r$ and
$s$. Hence, we have
$$\begin{aligned}
{\rm det}\begin{pmatrix} A & C \\ 0 & B \end{pmatrix}
&=(-1)^r(-1)^s{\rm det}\begin{pmatrix} A' & C' \\ 0 & B'
\end{pmatrix}
=(-1)^r(-1)^s{\rm det}\,A'\,\,{\rm det}\,B'\\
&=(-1)^r(-1)^s((-1)^r{\rm det}\,A)((-1)^s{\rm det}\,B)={\rm
det}\,A\,\,{\rm det}\,B
\end{aligned}$$ for some $C' \in M_{m \times n}(F)$. Also, we have
$${\rm det}\begin{pmatrix} A & 0 \\ D & B \end{pmatrix}
={\rm det}\begin{pmatrix} A & 0 \\ D & B
\end{pmatrix}^t={\rm det}\begin{pmatrix} A^t & D^t \\ 0 & B^t
\end{pmatrix}={\rm det}\,A^t\,\,{\rm det}\,B^t={\rm det}\,A\,\,{\rm det}\,B.$$
\end{proof}

Suppose that $A$ and $B$ are finite dimensional vector spaces over
a field $F$ and $f:A \rightarrow B$ is a linear transformation.
Then we have the short exact sequence
$$\begin{CD}
0 @>>> {\rm Ker}\,f @>\subseteq>> A @>f>> {\rm Im}\,f @>>> 0
\end{CD}\,.$$
If $f$ is 1-to-1, then ${\rm Ker}\,f=0$. Also, if $f$ is trivial,
then ${\rm Im}\,f=0$. Assume that $f$ is neither 1-to-1 nor
trivial. Let $k=(k_1,\dots,k_r)$, $b=(b_1,\dots,b_s)$, and
$a=(a_1,\dots,a_{r+s})$ be ordered bases for ${\rm Ker}\,f$, ${\rm
Im}\,f$, and $A$, respectively. Consider
$\widetilde{b}=(\widetilde{b}_1,\dots,\widetilde{b}_s)$ such that
$f(\widetilde{b}_i)=b_i$ for each $i \in \{1,\dots,s\}$. Then we
have an ordered basis
$(k,\widetilde{b})=(k_1,\dots,k_r,\widetilde{b}_1\dots,\widetilde{b}_s)$
for $A$, where
$\widetilde{b}=(\widetilde{b}_1,\dots,\widetilde{b}_s)$ is called
a lifting of $b$ by $f$, or simply, a lifting of $b$. If $b$ and
$b'$ are distinct ordered bases for ${\rm Im}\,f$, then
$\widetilde{b} \neq \widetilde{b'}$. However, although $b=b'$,
$\widetilde{b}$ and $\widetilde{b'}$ need not be the same. To
avoid this ambiguity, we write $\widetilde{b}^{(1)}$ for
$\widetilde{b}$ and $\widetilde{b}^{(2)}$ for $\widetilde{b'}$ if
$b=b'$.

\begin{lm} Let $f:A \rightarrow B$ be a linear transformation,
and let $k=(k_1,\dots,k_r)$ be an ordered basis for ${\rm
Ker}\,f$, and let $b=(b_1,\dots,b_s)$ be an ordered basis for
${\rm Im}\,f$. Then if $\widetilde{b}^{(1)}$ and
$\widetilde{b}^{(2)}$ are liftings of $b$, then
$$(k,\widetilde{b}^{(1)}) \sim (k,\widetilde{b}^{(2)}).$$
In other words, the equivalence class of $(k,\widetilde{b})$ with
respect to $\sim$ does not depend on the choice of lifting
$\widetilde{b}$ of $b$. We denote the equivalence class of
$(k,\widetilde{b})$ by $kb$.
\end{lm}

\begin{proof} The transition matrix $((k,\widetilde{b}^{(1)})/
(k,\widetilde{b}^{(2)}))$ is $\begin{pmatrix} I_r & 0 \\ C & I_s
\end{pmatrix}$ for some $C \in M_{s \times r}(F)$. Therefore,
$[(k,\widetilde{b}^{(1)})/(k,\widetilde{b}^{(2)})]=1$. That is,
$(k,\widetilde{b}^{(1)}) \sim (k,\widetilde{b}^{(2)})$.
\end{proof}

\begin{lm} Let $f:A \rightarrow B$ be a linear transformation,
and let $k=(k_1,\dots,k_r)$ and $k'=(k'_1,\dots,k'_r)$ be ordered
bases for ${\rm Ker}\,f$, and let $b=(b_1,\dots,b_s)$ and
$b'=(b'_1,\dots,b'_s)$ be distinct ordered bases for ${\rm
Im}\,f$. Then

(1) if $\widetilde{b}^{(1)}$ and $\widetilde{b}^{(2)}$ are
liftings of $b$, then
$[(k,\widetilde{b}^{(1)})/(k',\widetilde{b}^{(2)})]=[k/k']$;

(2) if $\widetilde{b}$ is a lifting of $b$ and $\widetilde{b'}$ is
a lifting of $b'$, then
$[(k,\widetilde{b})/(k,\widetilde{b'})]=[b/b']$.
\end{lm}

\begin{proof} (1) The transition matrix $((k,\widetilde{b}^{(1)})/
(k',\widetilde{b}^{(2)}))$ is $\begin{pmatrix} (k/k') & 0 \\ C &
I_s \end{pmatrix}$ for some $C \in M_{s \times r}(F)$. Therefore,
$[(k,\widetilde{b}^{(1)})/(k',\widetilde{b}^{(2)})]=[k/k']$. (2)
Similarly, the transition matrix $((k,\widetilde{b}^{(1)})/
(k,\widetilde{b}^{(2)}))$ is $\begin{pmatrix} I_r & 0 \\ C &
(b/b') \end{pmatrix}$ for some $C \in M_{s \times r}(F)$, so
$[(k,\widetilde{b})/(k,\widetilde{b'})]= [b/b']$.
\end{proof}

\begin{nota} By Lemma 2.2 and Lemma 2.4, we can express the
results of Lemma 2.5 in terms of equivalent classes as follows:

(1)
$[kb/k'b]=[(k,\widetilde{b}^{(1)})/(k',\widetilde{b}^{(2)})]=[k/k']$;
(2) $[kb/kb']=[(k,\widetilde{b})/(k,\widetilde{b'})]=[b/b']$.
\end{nota}

\begin{cor} Let $f:A \rightarrow B$ be a linear transformation.
Then if $k$ and $k'$ are ordered bases for ${\rm Ker}\,f$ and $b$
and $b'$ are distinct ordered bases for ${\rm Im}\,f$ and
$\widetilde{b}$ is a lifting of $b$ and $\widetilde{b'}$ is a
lifting of $b'$, then
$$[(k,\widetilde{b})/(k',\widetilde{b'})]=[k/k'][b/b'].$$
\end{cor}

\begin{proof} By Lemma 2.5, $[(k,\widetilde{b})/
(k',\widetilde{b'})]=[(k,\widetilde{b})/(k',\widetilde{b})]
[(k',\widetilde{b})/(k',\widetilde{b'})]=[k/k'][b/b']$.
\end{proof}

\begin{nota} The result of Corollary 2.7 can be written as
$$[kb/k'b']=[(k,\widetilde{b})/(k',\widetilde{b'})]=[k/k'][b/b'].$$
\end{nota}

\begin{lm} Let $A$ and $B$ be finite dimensional vector spaces
over a field $F$. Then if $a$ and $a'$ are ordered bases for $A$
and $b$ and $b'$ are ordered bases for $B$, then $(a,b)$ and
$(a',b')$ are ordered bases for $A \oplus B$ and
$[(a,b)/(a',b')]=[a/a'][b/b']$.
\end{lm}

\begin{proof} The transition matrix from $(a,b)$ to $(a',b')$ is
$((a,b)/(a',b'))=\begin{pmatrix} (a/a') & 0 \\ 0 & (b/b')
\end{pmatrix}$. Therefore, $[(a,b)/(a',b')]=[a/a'][b/b']$.
\end{proof}

We introduce the required definitions and properties from
algebraic topology and homological algebra. See \cite{AH}.

\begin{df} Let $C_0,\dots,C_m$ be modules over a ring $R$, and let
$\partial_i:C_{i+1} \rightarrow C_i$ be a $R$-module homomorphism
for each $i \in \{0,\dots,m-1\}$. Then
$$C=(\,\begin{CD}
0 @>>> C_m @>\partial_{m-1}>> C_{m-1} @>\partial_{m-2}>>
\cdots\cdots @>\partial_1>> C_1 @>\partial_0>> C_0 @>>> 0
\end{CD}\,)$$
is called a chain complex of length $m$ over $R$ if
$\partial_{i-1} \circ \partial_i=0$ for each $i \in
\{0,\dots,m\}$. Note that $C_{-1}=C_{m+1}=0$ and
$\partial_{-1}=\partial_m=0$. Also, we write the $R$-modules
$Z_i(C)={\rm Ker}\,\partial_{i-1}$, $B_i(C)={\rm Im}\,\partial_i$,
and $H_i(C)=Z_i(C)/B_i(C)$ for each $i \in \{0,\dots,m\}$ which
are called the $i$-th cycle, the $i$-th boundary, and the $i$-th
homology of the chain complex $C$, respectively. In particular, a
chain complex $C$ is said to be acyclic if $H_i(C)=0$ for each
$i$.
\end{df}

Note that $\partial_{i-1} \circ \partial_i=0$ if and only if ${\rm
Im}\,\partial_i \subseteq {\rm Ker}\,\partial_{i-1}$ and
$H_i(C)=0$ if and only if ${\rm Im}\,\partial_i={\rm
Ker}\,\partial_{i-1}$ for each $i \in \{0,\dots,m\}$.

\begin{df} A chain complex $C$ over a field $F$ is said to be
based if $C_i$ has a distinguished basis $c_i$ for each $i$.
\end{df}

Remark that we can think of a chain complex $C$ of length $m$ as
the chain complex $C$ of length $n$ for any $n \geq m$ by letting
$C_{m+1}=\dots=C_n=0$ and $\partial_m=\dots=\partial_{n-1}=0$. For
this reason, when we consider finitely many chain complexes with
different lengths simultaneously, we assume that they have the
same length $m$ which is the greatest length of them.

\begin{df} Let $C$ and $C'$ be chain complexes of length $m$
over a ring $R$. Then a sequence $f=(f_i:C_i \rightarrow
C'_i)_{i=0}^m$ of $R$-module homomorphisms is called a chain map
from $C$ to $C'$, denoted by $f:C \rightarrow C'$, if
$\partial'_{i-1} \circ f_i=f_{i-1} \circ \partial_{i-1}$ for each
$i \in \{1,\dots,m\}$. Note that $f_{-1}=0$ and $f_{m+1}=0$.
\end{df}

$$\begin{CD}
\cdots\cdots @>\partial_{i+1}>> C_{i+1} @>\partial_i>> C_i
@>\partial_{i-1}>> C_{i-1}
@>\partial_{i-2}>> \cdots\cdots\\
@VVV @Vf_{i+1}VV @Vf_iVV @Vf_{i-1}VV @VVV\\
\cdots\cdots @>\partial'_{i+1}>> C'_{i+1} @>\partial'_i>> C'_i
@>\partial'_{i-1}>> C'_{i-1}
@>\partial'_{i-2}>> \cdots\cdots\\
\end{CD}$$\\

\begin{pro} Let $C$ and $C'$ be chain complexes of length $m$
over a ring $R$. Then a chain map $f:C \rightarrow C'$ induces a
unique sequence $f_*=(f_{i*}:H_i(C) \rightarrow H_i(C'))_{i=0}^m$
of homomorphisms, denoted by $f_*:H_*(C) \rightarrow H_*(C')$,
such that for each $i \in \{0,\dots,m\}$, $f_{i*}:H_i(C)
\rightarrow H_i(C')$ is defined by $f_{i*}([z])=[f_i(z)]$ for all
$z \in Z_i(C)$. We call $f_*:H_*(C) \rightarrow H_*(C')$ the
induced homomorphism of $f$.
\end{pro}

\begin{proof} It suffices to show that for each $i \in \{0,\dots,m\}$,
$f_i(z) \in Z_i(C')$ for all $z \in Z_i(C)$ and $f_i(b) \in
B_i(C')$ for all $b \in B_i(C)$. Let $z \in Z_i(C)$. Then
$(f_{i-1} \circ \partial_{i-1})(z)=f_{i-1}(0)=0$, hence,
$(\partial'_{i-1} \circ f_i)(z)=0$. That is, $f_i(z) \in Z_i(C')$.
Let $b \in B_i(C)$. Then $b=\partial_i(c)$ for some $c \in
C_{i+1}$, hence, $f_i(b)=(f_i \circ \partial_i)(c)=(\partial'_i
\circ f_{i+1})(c)=\partial'_i(f_{i+1}(c))$. Since $f_{i+1}(c) \in
C'_{i+1}$, $f_i(b) \in B_i(C')$. Therefore, we have a unique
homomorphism $f_{i*}:H_i(C) \rightarrow H_i(C')$.
\end{proof}

\begin{df} Let $C$ and $C'$ be chain complexes of length $m$ over
a ring $R$, and let $f:C \rightarrow C'$ and $g:C \rightarrow C'$
be chain maps. Then $f$ and $g$ are said to be chain homotopic,
denoted by $f \simeq g$, if there is a sequence $T=(T_i:C_i
\rightarrow C'_{i+1})_{i=-1}^m$ of homomorphisms such that $f_i -
g_i=\partial'_i \circ T_i + T_{i-1} \circ
\partial_{i-1}$ for each $i \in \{0,\dots,m\}$. Note that $T_{-1}=0$.
When such a $T$ is known, which is called a chain homotopy between
$f$ and $g$, we say that $f$ and $g$ are chain homotopic by $T$.
\end{df}

\begin{pro} If $C$ and $C'$ are chain complexes of length $m$ over
a ring $R$, then the chain homotopic relation $\simeq$ on the set
$[C,C']$ of all chain maps from $C$ to $C'$ is an equivalence
relation.
\end{pro}

\begin{proof} For each $f \in [C,C']$, $f \simeq f$ by
zero map. That is, $T=0$. If $f \simeq g$ in $[C,C']$ by $T$, then
$g \simeq f$ by $-T$. If $f \simeq g$ in $[C,C']$ by $T_1$ and $g
\simeq h$ in $[C,C']$ by $T_2$, then $f \simeq h$ by $T_1 + T_2$.
Hence, $\simeq$ is an equivalence relation on $[C,C']$.
\end{proof}

\begin{pro} Let $C$, $C'$, and $C''$ be chain complexes of length $m$
over a ring $R$, and let $f:C \rightarrow C'$, $g:C \rightarrow
C'$, $f':C' \rightarrow C''$, and $g':C' \rightarrow C''$ be chain
maps. Then if $f \simeq g$ and $f' \simeq g'$, then $f' \circ f
\simeq g' \circ g:C \rightarrow C''$.
\end{pro}

\begin{proof} Suppose that $f \simeq g$ and
$T=(T_i:C_i \rightarrow C'_{i+1})_{i=-1}^m$ is a sequence of
homomorphisms such that $f_i - g_i=\partial'_i \circ T_i + T_{i-1}
\circ \partial_{i-1}$ for each $i \in \{0,\dots,m\}$. Let $i \in
\{0,\dots,m\}$ and $c \in C_i$. Then $f_i(c) -
g_i(c)=\partial'_i(T_i(c)) + T_{i-1}(\partial_{i-1}(c))$. Hence,
$$\begin{aligned}
(f'_i \circ f_i)(c) - (f'_i \circ g_i)(c)
&=((f'_i \circ \partial'_i) \circ T_i)(c) + ((f'_i \circ T_{i-1})
\circ \partial_{i-1})(c)\\
&=((\partial''_i \circ f'_{i+1}) \circ T_i)(c) + ((f'_i \circ
T_{i-1}) \circ \partial_{i-1})(c)\\
&=(\partial''_i \circ (f'_{i+1} \circ T_i))(c) + ((f'_i \circ
T_{i-1}) \circ \partial_{i-1})(c).
\end{aligned}$$
Therefore, $f' \circ f \simeq f' \circ g:C \rightarrow C''$ by $f'
\circ T = (f'_{i+1} \circ T_i:C_i \rightarrow
C''_{i+1})_{i=-1}^m$.

Similarly, we show that $f' \circ g \simeq g' \circ g:C
\rightarrow C''$.

Suppose that $f' \simeq g'$ and $T'=(T'_i:C'_i \rightarrow
C''_{i+1})_{i=-1}^m$ is a sequence of homomorphisms such that
$f'_i - g'_i=\partial''_i \circ T'_i + T'_{i-1} \circ
\partial'_{i-1}$ for each $i \in \{0,\dots,m\}$. Let $i \in
\{0,\dots,m\}$ and $c \in C_i$. Hence,
$$\begin{aligned}
(f'_i \circ g_i)(c) - (g'_i \circ g_i)(c)
&=(\partial''_i \circ T'_i)(g_i(c)) +
(T'_{i-1} \circ \partial'_{i-1})(g_i(c))\\
&=(\partial''_i \circ (T'_i \circ g_i))(c) + (T'_{i-1} \circ
(\partial'_{i-1} \circ g_i))(c)\\
&=(\partial''_i \circ (T'_i \circ g_i))(c) + (T'_{i-1} \circ
(g_{i-1} \circ \partial_{i-1}))(c)\\
&=(\partial''_i \circ (T'_i \circ g_i))(c) + ((T'_{i-1} \circ
g_{i-1}) \circ \partial_{i-1})(c).
\end{aligned}$$
Therefore, $f' \circ g \simeq g' \circ g:C \rightarrow C''$ by $T'
\circ g = (T'_i \circ g_i:C_i \rightarrow C''_{i+1})_{i=-1}^m$.
Furthermore, $f' \circ f \simeq g' \circ g:C \rightarrow C''$ by
$f' \circ T + T' \circ g$.
\end{proof}

\begin{pro} Let $C$ and $C'$ be chain complexes of length $m$ over
a field $F$. Then if chain maps $f:C \rightarrow C'$ and $g:C
\rightarrow C'$ are chain homotopic, then the induced
homomorphisms $f_* = g_*:H_*(C) \rightarrow H_*(C')$.
\end{pro}

\begin{proof} Suppose that $f$ and $g$ are chain homotopic and
$T=(T_i:C_i \rightarrow C'_{i+1})_{i=-1}^m$ is a sequence of
homomorphisms such that $f_i - g_i=\partial'_i \circ T_i + T_{i-1}
\circ \partial_{i-1}$ for each $i \in \{0,\dots,m\}$. Let $i \in
\{0,\dots,m\}$ and $z \in Z_i(C)$. Then $f_i(z) -
g_i(z)=\partial'_i(T_i(z)) +
T_{i-1}(\partial_{i-1}(z))=\partial'_i(T_i(z)) +
T_{i-1}(0)=\partial'_i(T_i(z)) \in B_i(C').$ Hence,
$f_{i*}([z])=[f_i(z)]=[g_i(z)]=g_{i*}([z])$. Therefore, $f_*=g_*$.
\end{proof}

\begin{df} Let $C$ and $C'$ be chain complexes of length $m$ over
a field $F$. Then $C$ and $C'$ are said to be chain equivalent if
there are chain maps $f:C \rightarrow C'$ and $g:C' \rightarrow C$
such that $g \circ f \simeq I_C$ and $f \circ g \simeq I_{C'}$,
where $I_C:C \rightarrow C$ and $I_{C'}:C' \rightarrow C'$ are the
identity chain maps. When such $f$ and $g$ are known, which are
called the chain equivalences between $C$ and $C'$, we say that
$C$ and $C'$ are chain equivalent by $(f,g)$.
\end{df}

\begin{pro} The chain equivalent relation $\simeq$ on the set
$K^m$ of all chain complexes of length $m$ over a ring $R$ is an
equivalence relation.
\end{pro}

\begin{proof} For each $C \in K^m$, $C \simeq C$ by
$(I_C,I_C)$. If $C \simeq C'$ in $K^m$ by $(f,g)$, then $C' \simeq
C$ by $(g,f)$. If $C \simeq C'$ in $K^m$ by $(f,g)$ and $C' \simeq
C''$ in $K^m$ by $(h,k)$, then $C \simeq C''$ by $(h \circ f,k
\circ g)$ by Proposition 2.16.
\end{proof}

\section{The torsion of a quasi-isomorphism}

In this section, we define a quasi-isomorphism and the torsion of
it.

\begin{df} Let $C$ and $C'$ be chain complexes of length $m$ over
a ring $R$. Then a chain map $f:C \rightarrow C'$ is said to be a
quasi-isomorphism if the induced homomorphism $f_*:H_*(C)
\rightarrow H_*(C')$ between homology is an isomorphism. That is,
the induced homomorphism $f_{i*}:H_i(C) \rightarrow H_i(C')$ is an
isomorphism for each $i \in \{0,\dots,m\}$.
\end{df}

\begin{pro} Let $C$ and $C'$ be chain complexes of length $m$ over
a ring $R$. Then if a chain map $f:C \rightarrow C'$ is an
isomorphism, then $f:C \rightarrow C'$ is a quasi-isomorphism.
\end{pro}

\begin{proof} We show that the induced homomorphism
$f_{i*}:H_i(C) \rightarrow H_i(C')$ is an isomorphism for each
$i$. Suppose that $z \in f_i^{-1}(B_i(C')) \cap Z_i(C)$. Then
$f_i(z) \in B_i(C')$ and $\partial_{i-1}(z)=0$. Let $w \in
D_{i+1}$ such that $f_i(z)=\partial'_i(w)$. Since $f_{i+1}$ is
onto, we can choose $x \in C_{i+1}$ so that $w=f_{i+1}(x)$. Hence,
$\partial'_i(f_{i+1}(x))=f_i(z)=f_i(\partial_i(x))$. Since $f_i$
is 1-to-1, $z=\partial_i(x) \in B_i(C)$. Hence, $f_{i*}$ is
1-to-1.

Suppose that $z' \in Z_i(C')$. Since $f_i$ is onto, we can choose
$z \in C_i$ so that $z'=f_i(z)$. Then
$\partial'_{i-1}(f_i(z))=0=f_{i-1}(\partial_{i-1}(z))$. Since
$f_{i-1}$ is 1-to-1, $\partial_{i-1}(z)=0$, that is, $z \in
Z_i(C)$. Hence $f_i:Z_i(C) \rightarrow Z_i(C')$ is onto.
Therefore, $f_{i*}$ is onto.
\end{proof}

To define the torsion of a quasi-isomorphism $f:C \rightarrow C'$,
we use the following short exact sequences
$$\begin{CD}
0 @>>> Z_i(C) @>\subseteq>> C_i @>\partial_{i-1}>> B_{i-1}(C) @>>>
0 \end{CD}\,,$$
$$\begin{CD}
0 @>>> B_i(C) @>\subseteq>> Z_i(C) @>\pi>> H_i(C) @>>> 0
\end{CD}\,,$$
$$\begin{CD}
0 @>>> Z_i(C') @>\subseteq>> C'_i @>\partial'_{i-1}>> B_{i-1}(C')
@>>> 0 \end{CD}\,,$$
$$\begin{CD}
0 @>>> B_i(C') @>\subseteq>> Z_i(C') @>\pi>> H_i(C') @>>> 0
\end{CD}$$\\
for each $i$, where $\pi$ is the canonical map.

\begin{df} Let $C$ and $C'$ be based chain complexes of length $m$
over a field $F$ such that $C_i$ and $C'_i$ are finite dimensional
vector spaces for each $i \in \{0,\dots,m\}$, and let $f:C
\rightarrow C'$ be a quasi-isomorphism. Then the torsion $\tau(f)$
of $f$ is defined by
$$\tau(f)=\prod_{i=0}^m
\left(\frac{[(b_ih_i)b_{i-1}/c_i]}{[(b'_if_{i*}(h_i))b'_{i-1}/c'_i]}
\right)^{(-1)^{i+1}},$$ where $b_i$, $b_{i-1}$, $b'_i$,
$b'_{i-1}$, $c_i$, $c'_i$, and $h_i$ are bases for $B_i(C)$,
$B_{i-1}(C)$, $B_i(C')$, $B_{i-1}(C')$, $C_i$, $C'_i$, and
$H_i(C)$, respectively, for each $i \in \{0,\dots,m\}$. Note that
$b_{-1}=b_m=b'_{-1}=b'_m=\emptyset$.
\end{df}

We can also write the torsion $\tau(f)$ of $f$ by
$$\tau(f)=\prod_{i=0}^m
\left(\frac{[(b_i,\,\widetilde{h_i},\,\widetilde{b_{i-1}})/c_i]}
{[(b'_i,\,f_i(\widetilde{h_i}),\,\widetilde{b'_{i-1}})/c'_i]}
\right)^{(-1)^{i+1}},$$ where $b_i$, $b_{i-1}$, $b'_i$,
$b'_{i-1}$, $c_i$, $c'_i$, and $h_i$ are bases for $B_i(C)$,
$B_{i-1}(C)$, $B_i(C')$, $B_{i-1}(C')$, $C_i$, $C'_i$, and
$H_i(C)$, respectively, and $\widetilde{b_{i-1}}$,
$\widetilde{b'_{i-1}}$, $\widetilde{h_i}$, and
$f_i(\widetilde{h_i})$ are liftings of $b_{i-1}$, $b'_{i-1}$,
$h_i$, and $f_{i*}(h_i)$, respectively, for each $i \in
\{0,\dots,m\}$.

$$\begin{CD}
0 @>>> C_m @>\partial_{m-1}>> C_{m-1} @>\partial_{m-2}>>
\cdots\cdots @>\partial_1>> C_1 @>\partial_0>> C_0 @>>> 0\\
@VVV @Vf_mVV @Vf_{m-1}VV @VVV @Vf_1VV @Vf_0VV @VVV\\
0 @>>> C'_m @>\partial'_{m-1}>> C'_{m-1} @>\partial'_{m-2}>>
\cdots\cdots @>\partial'_1>> C'_1 @>\partial'_0>> C'_0 @>>> 0
\end{CD}$$\\

\begin{lm} $\tau(f)$ dose not depend on the choices of $b_i$,
$b'_i$, and $h_i$. That is, the torsion $\tau$ on
quasi-isomorphisms is well-defined.
\end{lm}

\begin{proof} We use Notation 2.6 and 2.8 to prove this lemma.

For each $i \in \{0,\dots,m\}$, let $b_i$, $h_i$, and $b'_i$ be
bases for $B_i(C)$, $H_i(C)$, and $B_i(C')$, respectively. Then
$$\tau(f)=\prod_{i=0}^m
\left(\frac{[(b_ih_i)b_{i-1}/c_i]}{[(b'_if_{i*}(h_i))b'_{i-1}/c'_i]}
\right)^{(-1)^{i+1}}.$$

Step 1. Show that $\tau(f)$ is independent of the choice of $h_i$.

Let $i \in \{0,\dots,m\}$, and let $h'_i$ be a basis for $H_i(C)$.
Then we have

$[(b_ih'_i)b_{i-1}/c_i]=[(b_ih'_i)b_{i-1}/(b_ih_i)b_{i-1}]
[(b_ih_i)b_{i-1}/c_i]=[b_ih'_i/b_ih_i] [(b_ih_i)b_{i-1}/c_i]$

$=[h'_i/h_i][(b_ih_i)b_{i-1}/c_i]$ and

$[(b'_if_{i*}(h'_i))b'_{i-1}/c'_i]
=[(b'_if_{i*}(h'_i))b'_{i-1}/(b'_if_{i*}(h_i))b'_{i-1}]
[(b'_if_{i*}(h_i))b'_{i-1}/c'_i]$

$=[b'_if_{i*}(h'_i)/b'_if_{i*}(h_i)][(b'_if_{i*}(h_i))b'_{i-1}/c'_i]$

$=[f_{i*}(h'_i)/f_{i*}(h_i)][(b'_if_{i*}(h_i))b'_{i-1}/c'_i]$.

Since $f_{i*}$ is an isomorphism,
$(h'_i/h_i)=(f_{i*}(h'_i)/f_{i*}(h_i))$, so
$[h'_i/h_i]=[f_{i*}(h'_i)/f_{i*}(h_i)]$. Hence, we have
$$\frac{[(b_ih'_i)b_{i-1}/c_i]}{[(b'_if_{i*}(h'_i))b'_{i-1}/c'_i]}
=\frac{[(b_ih_i)b_{i-1}/c_i]}{[(b'_if_{i*}(h_i))b'_{i-1}/c'_i]}$$
for each $i \in \{0,\dots,m\}$.

Therefore, $\tau(f)$ does not depend on the choice of bases for
homology spaces.

Step 2. Show that $\tau(f)$ is independent of the choices of $b_i$
and $b'_i$.

For each $i \in \{0,\dots,m-1\}$, let
$X_i(b_{i-1},b'_{i-1},b_i,b'_i,b_{i+1},b'_{i+1})$ be the
expression
$$\left(\frac{[(b_ih_i)b_{i-1}/c_i]}{[(b'_if_{i*}(h_i))b'_{i-1}/c'_i]}
\right)^{(-1)^{i+1}}
\left(\frac{[(b_{i+1}h_{i+1})b_i/c_{i+1}]}{[(b'_{i+1}f_{i+1*}(h_{i+1}))
b'_i/c'_{i+1}]} \right)^{(-1)^{i+2}}.$$ We claim that
$X_i(b_{i-1},b'_{i-1},b_i,b'_i,b_{i+1},b'_{i+1})$ is independent
of the choices of $b_i$ and $b'_i$. Let $i \in \{0,\dots,m-1\}$,
and let $v_i$ and $v'_i$ be bases for $B_i(C)$ and $B_i(C')$,
respectively. Then we have the following equations

$[(v_ih_i)b_{i-1}/c_i]=[(v_ih_i)b_{i-1}/(b_ih_i)b_{i-1}]
[(b_ih_i)b_{i-1}/c_i]=[v_ih_i/b_ih_i][(b_ih_i)b_{i-1}/c_i]$

$=[v_i/b_i][(b_ih_i)b_{i-1}/c_i]$,\\

$[(b_{i+1}h_{i+1})v_i/c_{i+1}]=[(b_{i+1}h_{i+1})v_i/(b_{i+1}h_{i+1})b_i]
[(b_{i+1}h_{i+1})b_i/c_{i+1}]$

$=[b_{i+1}h_{i+1}/b_{i+1}h_{i+1}][v_i/b_i][(b_{i+1}h_{i+1})b_i/c_{i+1}]$

$=[v_i/b_i][(b_{i+1}h_{i+1})b_i/c_{i+1}]$,\\

$[(v'_if_{i*}(h_i))b'_{i-1}/c'_i]
=[(v'_if_{i*}(h_i))b'_{i-1}/(b'_if_{i*}(h_i))b'_{i-1}]
[(b'_if_{i*}(h_i))b'_{i-1}/c'_i]$

$=[v'_if_{i*}(h_i)/b'_if_{i*}(h_i)][(b'_if_{i*}(h_i))b'_{i-1}/c'_i]$

$=[v'_i/b'_i][(b'_if_{i*}(h_i))b'_{i-1}/c'_i]$,\\

$[(b'_{i+1}f_{i+1*}(h_{i+1}))v'_i/c'_{i+1}]$

$=[(b'_{i+1}f_{i+1*}(h_{i+1}))v'_i/(b'_{i+1}f_{i+1*}(h_{i+1}))b'_i]
[(b'_{i+1}f_{i+1*}(h_{i+1}))b'_i/c'_{i+1}]$

$=[b'_{i+1}f_{i+1*}(h_{i+1})/b'_{i+1}f_{i+1*}(h_{i+1})][v'_i/b'_i]
[(b'_{i+1}f_{i+1*}(h_{i+1}))b'_i/c'_{i+1}]$

$=[v'_i/b'_i][(b'_{i+1}f_{i+1*}(h_{i+1}))b'_i/c'_{i+1}]$.\\

Hence, we have
$$X_i(b_{i-1},b'_{i-1},v_i,v'_i,b_{i+1},b'_{i+1})=
X_i(b_{i-1},b'_{i-1},b_i,b'_i,b_{i+1},b'_{i+1})$$ for each $i \in
\{0,\dots,m-1\}$. Also, we have
$$X_i(v_{i-1},v'_{i-1},v_i,v'_i,b_{i+1},b'_{i+1})=
X_i(v_{i-1},v'_{i-1},b_i,b'_i,b_{i+1},b'_{i+1})$$ for each $i \in
\{0,\dots,m-1\}$. Remark that $v_{-1}=v_m=v'_{-1}=v'_m=\emptyset$.
Therefore, by these facts, we conclude that $$\prod_{i=0}^m
\left(\frac{[(v_ih_i)v_{i-1}/c_i]}{[(v'_if_{i*}(h_i))v'_{i-1}/c'_i]}
\right)^{(-1)^{i+1}}=\prod_{i=0}^m
\left(\frac{[(b_ih_i)b_{i-1}/c_i]}{[(b'_if_{i*}(h_i))b'_{i-1}/c'_i]}
\right)^{(-1)^{i+1}}.$$

This proves the lemma.
\end{proof}

\begin{lm} Let $C$, $C'$, and $C''$ be based chain complexes of
length $m$ over a field $F$. Then if $f:C \rightarrow C'$ and
$g:C' \rightarrow C''$ are quasi-isomorphisms, then $g \circ f:C
\rightarrow C''$ is a quasi-isomorphism and $\tau(g \circ
f)=\tau(g)\tau(f)$.
\end{lm}

\begin{proof}
$$\begin{CD}
0 @>>> C_m @>\partial_{m-1}>> C_{m-1} @>\partial_{m-2}>>
\cdots\cdots @>\partial_1>> C_1 @>\partial_0>> C_0 @>>> 0\\
@VVV @Vf_mVV @Vf_{m-1}VV @VVV @Vf_1VV @Vf_0VV @VVV\\
0 @>>> C'_m @>\partial'_{m-1}>> C'_{m-1} @>\partial'_{m-2}>>
\cdots\cdots @>\partial'_1>> C'_1 @>\partial'_0>> C'_0 @>>> 0\\
@VVV @Vg_mVV @Vg_{m-1}VV @VVV @Vg_1VV @Vg_0VV @VVV\\
0 @>>> C''_m @>\partial''_{m-1}>> C''_{m-1} @>\partial''_{m-2}>>
\cdots\cdots @>\partial''_1>> C''_1 @>\partial''_0>> C''_0 @>>> 0
\end{CD}$$\\
For each $i \in \{0,\dots,m\}$, $\partial''_{i-1} \circ (g_i \circ
f_i)=(g_{i-1} \circ f_{i-1}) \circ \partial_{i-1}$ and $(g_i \circ
f_i)_*=g_{i*} \circ f_{i*}:H_i(C) \rightarrow H_i(C'')$ is an
isomorphism. Hence, $g \circ f:C \rightarrow C''$ is a
quasi-isomorphism. Note that the torsion does not depend on the
choice of basis for homology.

Suppose that $$\tau(f)=\prod_{i=0}^m
\left(\frac{[(b_ih_i)b_{i-1}/c_i]}{[(b'_if_{i*}(h_i))b'_{i-1}/c'_i]}
\right)^{(-1)^{i+1}}.$$ Then $$\tau(g)=\prod_{i=0}^m
\left(\frac{[(b'_if_{i*}(h_i))b'_{i-1}/c'_i]}
{[(b''_ig_{i*}(f_{i*}(h_i)))b''_{i-1}/c''_i]}
\right)^{(-1)^{i+1}}.$$ Hence,
$$\begin{aligned}
\tau(g)\tau(f) &=\tau(f)\tau(g)=\prod_{i=0}^m
\left(\frac{[(b_ih_i)b_{i-1}/c_i]}
{[(b''_ig_{i*}(f_{i*}(h_i)))b''_{i-1}/c''_i]}
\right)^{(-1)^{i+1}}\\
&=\prod_{i=0}^m \left(\frac{[(b_ih_i)b_{i-1}/c_i]} {[(b''_i(g_i
\circ f_i)_*(h_i))b''_{i-1}/c''_i]} \right)^{(-1)^{i+1}}=\tau(g
\circ f).
\end{aligned}$$
\end{proof}

\begin{lm} Let $C$, $C'$, $C''$, and $C'''$ be based chain complexes
of length $m$ over a field $F$. Then if $f:C \rightarrow C'$ and
$g:C'' \rightarrow C'''$ are quasi-isomorphisms, then $f \oplus
g:C \oplus C'' \rightarrow C' \oplus C'''$ is a quasi-isomorphism
and $\tau(f \oplus g)=\pm\,\tau(f)\tau(g)$.
\end{lm}

\begin{proof} For each $i \in \{0,\dots,m\}$, $(\partial'_{i-1} \oplus
\partial'''_{i-1}) \circ (f_i \oplus g_i)=(\partial'_{i-1} \circ f_i)
\oplus (\partial'''_{i-1} \circ g_i)=(f_{i-1} \circ
\partial_{i-1}) \oplus (g_{i-1} \circ \partial''_{i-1})=(f_{i-1}
\oplus g_{i-1}) \circ (\partial_{i-1} \oplus \partial''_{i-1})$
and $(f_i \oplus g_i)_*: H_i(C \oplus C'') \rightarrow H_i(C'
\oplus C''')$ is an isomorphism since $f_{i*} \oplus g_{i*}:H_i(C)
\oplus H_i(C'') \rightarrow H_i(C') \oplus H_i(C''')$ is an
isomorphism. Hence, $f \oplus g:C \oplus C'' \rightarrow C' \oplus
C'''$ is a quasi-isomorphism. Also,
$$\begin{aligned}
\tau(f \oplus g) &=\prod_{i=0}^m \left(\frac{[((b_i \oplus
b''_i)(h_i \oplus h''_i))(b_{i-1} \oplus b''_{i-1})/(c_i \oplus
c''_i)]}{[((b'_i \oplus b'''_i)(f_i \oplus g_i)_*(h_i \oplus
h''_i))(b'_{i-1} \oplus b'''_{i-1})/(c'_i
\oplus c'''_i)]} \right)^{(-1)^{i+1}}\\
&=\prod_{i=0}^m \left(\frac{[((b_i \oplus b''_i)(h_i \oplus
h''_i))(b_{i-1} \oplus b''_{i-1})/(c_i \oplus c''_i)]}{[((b'_i
\oplus b'''_i)(f_{i*}(h_i) \oplus g_{i*}(h''))(b'_{i-1} \oplus
b'''_{i-1})/(c'_i \oplus c'''_i)]} \right)^{(-1)^{i+1}}\\
&=\pm\,\prod_{i=0}^m \left(\frac{[((b_ih_i)b_{i-1} \oplus
(b''_ih''_i)b''_{i-1})/(c_i \oplus
c''_i)]}{[(b'_if_{i*}(h_i))b'_{i-1} \oplus
(b'''_ig_{i*}(h''_i))b'''_{i-1}/(c'_i \oplus c'''_i)]}
\right)^{(-1)^{i+1}}\\
&=\pm\,\prod_{i=0}^m
\left(\frac{[(b_ih_i)b_{i-1}/c_i]}{[(b'_if_{i*}(h_i))b'_{i-1}/c'_i]}
\frac{[(b''_ih''_i)b''_{i-1}/c''_i]}
{[(b'''_ig_{i*}(h''_i))b'''_{i-1}/c'''_i]}
\right)^{(-1)^{i+1}}\\
&=\pm\,\prod_{i=0}^m
\left(\frac{[(b_ih_i)b_{i-1}/c_i]}{[(b'_if_{i*}(h_i))b'_{i-1}/c'_i]}
\right)^{(-1)^{i+1}} \prod_{i=0}^m
\left(\frac{[(b''_ih''_i)b''_{i-1}/c''_i]}
{[(b'''_ig_{i*}(h''_i))b'''_{i-1}/c'''_i]}
\right)^{(-1)^{i+1}}\\
&=\pm\,\tau(f)\tau(g).
\end{aligned}$$
\end{proof}

Notice that a sign problem occurs at the 3rd equation. The
dimensions of $C_i$ and $C'_i$ are not the same, even boundaries
$B_i(C)$ and $B_i(C')$, in general.

To get the exact sign for the torsion of it, for each $i \in
\{0,\dots,m\}$, let $$x_i={\rm dim}_F\,B_i(C), \,\,\,x'_i={\rm
dim}_F\,B_i(C'), \,\,\,x''_i={\rm dim}_F\,B_i(C''),
\,\,\,x'''_i={\rm dim}_F\,B_i(C'''),$$
$$y_i={\rm dim}_F\,H_i(C), \,\,\,y''_i={\rm dim}_F\,H_i(C''),$$ Then
$$\tau(f \oplus g)=\frac{(-1)^{\sum_{i=0}^m(x''_iy_i + x_{i-1}
(x''_i + y''_i))}}{(-1)^{\sum_{i=0}^m(x'''_iy_i + x'_{i-1}(x'''_i
+ y''_i))}}\,\tau(f)\tau(g).$$ Therefore,
$$\tau(f \oplus g)=(-1)^{\sum_{i=0}^m[(x''_iy_i + x_{i-1}
(x''_i + y''_i))-(x'''_iy_i + x'_{i-1}(x'''_i +
y''_i))]}\,\tau(f)\tau(g).$$

In particular, if $C=C'$ and $C''=C'''$, then $\tau(f \oplus
g)=\tau(f)\tau(g)$.

\begin{cor} Let $C$, $C'$, $C''$, and $C'''$ be based chain complexes
of length $m$ over a field $F$. Then if $f:C \rightarrow C'$ and
$g:C'' \rightarrow C'''$ are quasi-isomorphisms, then $\tau(g
\oplus f)=\pm\,\tau(f \oplus g)$.
\end{cor}

\begin{proof} Since $\tau(f \oplus g)=\pm\,\tau(f)\tau(g)$ and
$\tau(g \oplus f)=\pm\,\tau(g)\tau(f)$, we have $\tau(g \oplus
f)=\pm\,\tau(f \oplus g)$.
\end{proof}

\begin{cor} Let $C$, $C'$, $C''$, $C'''$, $C''''$, and $C'''''$ be
based chain complexes of length $m$ over a field $F$. Then if $f:C
\rightarrow C'$, $g:C'' \rightarrow C'''$, and $h:C''''
\rightarrow C'''''$ are quasi-isomorphisms, then $\tau((f \oplus
g) \oplus h)=\pm\,\tau(f \oplus (g \oplus h))$.
\end{cor}

\begin{proof} Since $\tau((f \oplus g) \oplus h)=\pm\,\tau(f \oplus g)
\tau(h)=\pm\,\tau(f)\tau(g)\tau(h)$ and $\tau(f \oplus (g \oplus
h))=\pm\,\tau(f)\tau(g \oplus h)= \pm\,\tau(f)\tau(g)\tau(h)$, we
have $\tau((f \oplus g) \oplus h)=\pm\,\tau(f \oplus (g \oplus
h))$.
\end{proof}

\begin{lm} Let $C$ and $C'$ be based chain complexes of length $m$
over a field $F$. Then if quasi-isomorphisms $f:C \rightarrow C'$
and $g:C \rightarrow C'$ are chain homotopic, then
$\tau(f)=\tau(g)$.
\end{lm}

\begin{proof} Since $f$ and $g$ are chain homotopic, $f_*=g_*$.
Therefore, $\tau(f)=\tau(g)$.
\end{proof}

\begin{lm} Let $C$ and $C'$ be based chain complexes of length $m$
over a field $F$. Then if $f:C \rightarrow C'$ and $g:C'
\rightarrow C$ are the chain equivalences between $C$ and $C'$,
then $f$ and $g$ are quasi-isomorphisms and
$\tau(f)=\tau(g)^{-1}$.
\end{lm}

\begin{proof} Since $g \circ f \simeq I_{C}$ and $f \circ g \simeq
I_{C'}$, we have $g_* \circ f_* = (g \circ f)_* = I_{C*}$ and $f_*
\circ g_* = (f \circ g)_* = I_{C'*}$. Note that $I_{C*}$ and
$I_{C'*}$ are the identity induced maps. Hence, $f_*$ and $g_*$
are isomorphisms, so $f$ and $g$ are quasi-isomorphisms.
Therefore, $\tau(g)\tau(f)=\tau(g \circ f)=\tau(I_C)=1$. Hence, we
have $\tau(f)=\tau(g)^{-1}$.
\end{proof}

\begin{df} Let $C$ be a based chain complex of length $m$ over
a field $F$. Then chain maps $f:C \rightarrow C$ and $g:C
\rightarrow C$ are said to be conjugate if there is a chain
isomorphism $h:C \rightarrow C$ such that $f = h^{-1} \circ g
\circ h$.
\end{df}

\begin{lm} Let $C$ be a based chain complex of length $m$ over
a field $F$. Then if quasi-isomorphisms $f:C \rightarrow C$
and $g:C \rightarrow C$ are conjugate, then $\tau(f)=\tau(g)$.
\end{lm}

\begin{proof} Suppose that $f = h^{-1} \circ g \circ h$ for some
chain isomorphism $h:C \rightarrow C$. Then $\tau(f)=\tau(h^{-1}
\circ g \circ h)=\tau(h^{-1})\tau(g)\tau(h)$. Since $h$ is a chain
isomorphism, $h$ and $h^{-1}$ are chain equivalences. Hence, by
Lemma 3.10, we have $\tau(f)=\tau(g)$.
\end{proof}

Let us introduce the definition of torsion of a based acyclic
chain complex which gives us a motivation to define our torsion of
a quasi-isomorphism. See \cite{T}. The torsion of a based acyclic
chain complex can be regarded as the torsion of a
quasi-isomorphism.

\begin{df} Let $C$ be a based acyclic chain complex of length $m$
over a field $F$ such that $C_i$ is a finite dimensional vector
space for each $i \in \{0,\dots,m\}$. Then the torsion $\tau(C)$
of $C$ is defined by $$\tau(C)=\prod_{i=0}^m
[b_ib_{i-1}/c_i]^{(-1)^{i+1}},$$ where $b_i$ and $c_i$ are bases
for $B_i(C)$ and $C_i$, respectively, for each $i \in
\{0,\dots,m\}$. Note that $b_{-1}=b_m=\emptyset$. In particular,
the zero chain complex of length $m$, denoted by $0^m$, is acyclic
and we define $\tau(0^m)=1$.
\end{df}

\begin{thm} If $C$ and $C'$ are acyclic based chain complexes of
length $m$ over a field $F$ and $f:C \rightarrow C'$ is a chain
map, then $f$ is a quasi-isomorphism and
$$\tau(f)=\frac{\tau(C)}{\tau(C')}\,.$$
\end{thm}

\begin{proof} Since $C$ and $C'$ are acyclic, $H_i(C)=H_i(C')=0$
for each $i \in \{0,\dots,m\}$. Therefore, $f$ is a
quasi-isomorphism and $$\tau(f)=\prod_{i=0}^m
\left(\frac{[(b_ih_i)b_{i-1}/c_i]}{[(b'_if_{i*}(h_i))b'_{i-1}/c'_i]}
\right)^{(-1)^{i+1}}=\prod_{i=0}^m
\left(\frac{[b_ib_{i-1}/c_i]}{[b'_ib'_{i-1}/c'_i]}
\right)^{(-1)^{i+1}}=\frac{\tau(C)}{\tau(C')}\,,$$ where $b_i$,
$b_{i-1}$, $b'_i$, $b'_{i-1}$, $c_i$, and $c'_i$ are bases for
$B_i(C)$, $B_{i-1}(C)$, $B_i(C')$, $B_{i-1}(C')$, $C_i$, and
$C'_i$, respectively, and $h_i=\emptyset$ for each $i \in
\{0,\dots,m\}$.
\end{proof}

\begin{nota} Let $C$ be a based chain complex of length $m$ over
a field $F$. Then we denote zero chain maps which are injective
and surjective by
$$0_{0 \rightarrow C}:0^m \rightarrow C \,\,\,{\rm and}\,\,\,
0_{C \rightarrow 0}:C \rightarrow 0^m,$$ respectively.
\end{nota}

\begin{cor} If $C$ is a based acyclic chain complex of length $m$
over a field $F$, then $0_{0 \rightarrow C}:0^m \rightarrow C$ and
$0_{C \rightarrow 0}:C \rightarrow 0^m$ are quasi-isomorphisms and
$$\tau(C)=\tau(0_{0 \rightarrow C})^{-1}=\tau(0_{C \rightarrow 0}).$$
\end{cor}

\begin{proof} Since $C$ and $0_m$ are acyclic, by Theorem 3.14,
we have $$\tau(0_{0 \rightarrow
C})=\frac{\tau(0^m)}{\tau(C)}=\tau(C)^{-1} \,\,\,{\rm and}\,\,\,
\tau(0_{C \rightarrow 0})=\frac{\tau(C)}{\tau(0^m)}=\tau(C).$$
\end{proof}

\begin{thm} Let $C$ and $C'$ be based chain complexes of length $m$
over a field $F$. Then if $C'$ is acyclic, then the sequence $i:C
\rightarrow C \oplus C'$ of injection maps and the sequence $p:C
\oplus C' \rightarrow C$ of projection maps are quasi-isomorphisms
and $\tau(i)=\pm\,\tau(C')^{-1}$ and $\tau(p)=\pm\,\tau(C')$.
\end{thm}

\begin{proof} First, we show that $p$ and $i$ are
quasi-isomorphisms. Let $j \in \{0,\dots,m\}$ and $c \in C_j$ and
$c' \in C'_j$. Then $((\partial_{j-1} \oplus \partial'_{j-1})
\circ i_j)(c)=(i_{j-1} \circ \partial_{j-1})(c)$ and
$(\partial_{j-1} \circ p_j)(c \oplus c')=(p_{j-1} \circ
(\partial_{j-1} \oplus \partial'_{j-1}))(c \oplus c')$. Hence,
$i:C \rightarrow C \oplus C'$ and $p:C \oplus C' \rightarrow C$
are chain maps, which are called the injection chain map and the
projection chain map, respectively. We can think of $i:C
\rightarrow C \oplus C'$ and $p:C \oplus C' \rightarrow C$ as $I_C
\oplus 0_{0 \rightarrow C'}:C \oplus 0 \rightarrow C \oplus C'$
and $I_C \oplus 0_{C' \rightarrow 0}:C \oplus C' \rightarrow C
\oplus 0$, respectively. By Lemma 3.6,
$\tau(i)=\pm\,\tau(C')^{-1}$ and $\tau(p)=\pm\,\tau(C')$.
\end{proof}

We need to distinguish a chain complex with distinct bases. Let us
use a pair $(C,c)$ for a based chain complex $C$ with a basis $c$.
Also, a chain map $f:C \rightarrow C'$ between chain complexes $C$
and $C'$ with bases $c$ and $c'$, respectively, is denoted by
$f:(C,c) \rightarrow (C',c')$.

\begin{lm} Let $C$ be a chain complex of length $m$ over
a field $F$, and let $I_C:C \rightarrow C$ be the identity chain
map. Then if $c^{(1)}$ and $c^{(2)}$ are bases for $C$, then
$$\tau(I_C:(C,c^{(1)}) \rightarrow (C,c^{(2)}))=\prod_{i=0}^m
[c^{(2)}_i/c^{(1)}_i]^{(-1)^{i+1}}.$$
\end{lm}

\begin{proof} Suppose that $b_i$, $b_{i-1}$, and $h_i$ are bases
for $B_i(C)$, $B_{i-1}(C)$, and $H_i(C)$, respectively, for each
$i \in \{0,\dots,m\}$. Since the torsion is independent of the
choice of bases for boundaries, we have
$$\begin{aligned}
&\tau(I_C:(C,c^{(1)}) \rightarrow (C,c^{(2)}))=\prod_{i=0}^m
\left(\frac{[(b_ih_i)b_{i-1}/c^{(1)}_i]}{[(b_iI_{Ci*}(h_i))b_{i-1}/c^{(2)}_i]}
\right)^{(-1)^{i+1}}\\
&=\prod_{i=0}^m
\left(\frac{[(b_ih_i)b_{i-1}/c^{(1)}_i]}{[(b_ih_i)b_{i-1}/c^{(2)}_i]}
\right)^{(-1)^{i+1}}=\prod_{i=0}^m
[c^{(2)}_i/c^{(1)}_i]^{(-1)^{i+1}}.
\end{aligned}$$
\end{proof}

\begin{thm} Let $C$ and $C'$ be chain complexes of length $m$ over
a field $F$, and let $f:C \rightarrow C'$ be a quasi-isomorphism.
Then if $c^{(1)}$ and $c^{(2)}$ are bases for $C$ and $c'^{(1)}$
and $c'^{(2)}$ are bases for $C'$, then
$$\tau(f:(C,c^{(2)}) \rightarrow (C',c'^{(2)}))=\prod_{i=0}^m
\left(\frac{[c^{(1)}_i/c^{(2)}_i]}{[c'^{(1)}_i/c'^{(2)}_i]}\right)
^{(-1)^{i+1}} \tau(f:(C,c^{(1)}) \rightarrow (C',c'^{(1)})).$$
\end{thm}

\begin{proof} Suppose that $f^{(1)}=f:(C,c^{(1)}) \rightarrow
(C',c'^{(1)})$, $f^{(2)}=f:(C,c^{(2)}) \rightarrow (C',c'^{(2)})$,
$I^{(21)}=I_C:(C,c^{(2)}) \rightarrow (C,c^{(1)})$, and
$I'^{(12)}=I_{C'}:(C',c'^{(1)}) \rightarrow (C',c'^{(2)})$. Since
$f^{(2)}=I'^{(12)} \circ f^{(1)} \circ I^{(21)}$, we have
$\tau(f^{(2)})=\tau(I'^{(12)})\tau(f^{(1)})\tau(I^{(21)})$ by
Lemma 3.5. Hence, by Lemma 3.18,
$$\tau(f^{(2)})=\prod_{i=0}^m
[c'^{(2)}_i/c'^{(1)}_i]^{(-1)^{i+1}} \tau(f^{(1)}) \prod_{i=0}^m
[c^{(1)}_i/c^{(2)}_i]^{(-1)^{i+1}}.$$ Therefore,
$$\tau(f^{(2)})=\prod_{i=0}^m
\left(\frac{[c^{(1)}_i/c^{(2)}_i]}{[c'^{(1)}_i/c'^{(2)}_i]}\right)
^{(-1)^{i+1}} \tau(f^{(1)}).$$
\end{proof}

The torsion of a quasi-isomorphism from a based chain complex to
itself can be easily calculated. It turns out that the torsion is
determined by the determinants of the induced isomorphisms on
homology.

\begin{thm} Let $C$ be a based chain complex of length $m$ over
a field $F$, and let $f:C \rightarrow C$ be a quasi-isomorphism.
Then

(1) $\tau(f)$ is independent of the choice of the basis $c$;

(2) if $h_i$ is a basis for $H_i(C)$ for each $i \in
\{0,\dots,m\}$, then $$\tau(f)=\prod_{i=0}^m
[h_i/f_{i*}(h_i)]^{(-1)^{i+1}}=\frac{\prod_{i\,even}{{\rm
det}\,f_{i*}}}{\prod_{i\,odd} {\rm det}\,f_{i*}},$$ where ${\rm
det}\,f_{i*}=[f_{i*}(h_i)/h_i]$ for each $i \in \{0,\dots,m\}$.
\end{thm}

\begin{proof} (1) If $c^{(1)}$ and $c^{(2)}$ are bases for $C$,
then, by Theorem 3.19, $$\tau(f:(C,c^{(2)}) \rightarrow
(C,c^{(2)}))=\prod_{i=0}^m
\left(\frac{[c^{(1)}_i/c^{(2)}_i]}{[c^{(1)}_i/c^{(2)}_i]}\right)
^{(-1)^{i+1}} \tau(f:(C,c^{(1)}) \rightarrow (C,c^{(1)})).$$
Hence, $\tau(f:(C,c^{(2)}) \rightarrow
(C,c^{(2)}))=\tau(f:(C,c^{(1)}) \rightarrow (C,c^{(1)}))$.

(2) $$\begin{aligned} \tau(f) &=\prod_{i=0}^m
\left(\frac{[(b_ih_i)b_{i-1}/c_i]}{[(b_if_{i*}(h_i))b_{i-1}/c_i]}
\right)^{(-1)^{i+1}}=\prod_{i=0}^m
\left([(b_ih_i)b_{i-1}/c_i][c_i/(b_if_{i*}(h_i))b_{i-1}]
\right)^{(-1)^{i+1}}\\
&=\prod_{i=0}^m
[(b_ih_i)b_{i-1}/(b_if_{i*}(h_i))b_{i-1}]^{(-1)^{i+1}}
=\prod_{i=0}^m [h_i/f_{i*}(h_i)]^{(-1)^{i+1}}.
\end{aligned}$$
Also, if $i \in \{1,\dots,m\}$ is even, then
$[h_i/f_{i*}(h_i)]^{(-1)^{i+1}}=[f_{i*}(h_i)/h_i]$. Similarly, if
$i \in \{1,\dots,m\}$ is odd, then
$[h_i/f_{i*}(h_i)]^{(-1)^{i+1}}=[h_i/f_{i*}(h_i)]=[f_{i*}(h_i)/h_i]^{-1}$.
Notice that $[h_i/f_{i*}(h_i)]=({\rm det}\,f_{i*})^{-1}$ for each
$i \in \{0,\dots,m\}$. Therefore, $$\tau(f)=\prod_{i=0}^m
[h_i/f_{i*}(h_i)]^{(-1)^{i+1}}=\prod_{i=0}^m \left(\frac{1}{{\rm
det}\,f_{i*}} \right)^{(-1)^{i+1}}=\frac{\prod_{i\,even}{{\rm
det}\,f_{i*}}}{\prod_{i\,odd} {\rm det}\,f_{i*}}.$$
\end{proof}

The following theorem is a statement which can be regarded as the
generalization of Proposition 2.3 to the torsion of
quasi-isomorphisms.

\begin{thm} Let $C$ and $C'$ be based chain complexes of length $m$
over a field $F$, and let $f:C \rightarrow C$ and $f':C'
\rightarrow C'$ be quasi-isomorphisms, and let $g:C \rightarrow
C'$ be a chain map. Define $\bar{f}:C \oplus C' \rightarrow C
\oplus C'$ by $\bar{f}_i(x \oplus y)=f_i(x) \oplus
(g_i(x)+f'_i(y))$ for all $x \in C_i$ and $y \in C'_i$ for each $i
\in \{0,\dots,m\}$. Then $\bar{f}:C \oplus C' \rightarrow C \oplus
C'$ is a quasi-isomorphism and $\tau(\bar{f})=\tau(f)\tau(f')$. In
particular, if $g=0$, then $\bar{f}=f \oplus f'$.
\end{thm}

\begin{proof} For convenience, let us use $x \oplus y$ for an ordered
pair $(x,y)$.

Step 1. Show that $\bar{f}:C \oplus C' \rightarrow C \oplus C'$ is
a chain map.

Let $i \in \{0,\dots,m\}$. Suppose that $x,x_1,x_2 \in C$ and
$y,y_1,y_2 \in C'$ and $r \in F$. Then we have
$$\begin{aligned}
&\bar{f}_i(x_1 \oplus y_1 + x_2 \oplus y_2)=\bar{f}_i((x_1+x_2)
\oplus (y_1+y_2))\\
&=f_i(x_1+x_2) \oplus (g_i(x_1+x_2)+f'_i(y_1+y_2))\\
&=(f_i(x_1) \oplus (g_i(x_1)+f'_i(y_1))) + (f_i(x_2) \oplus
(g_i(x_2)+f'_i(y_2)))\\
&=\bar{f}_i(x_1 \oplus y_1) + \bar{f}_i(x_2 \oplus y_2)
\end{aligned}$$ and
$$\bar{f}_i(r(x \oplus y))=\bar{f}_i(rx \oplus ry)=f_i(rx) \oplus
(g_i(rx)+f'_i(ry))=r\bar{f}_i(x \oplus y)$$ and
$$\begin{aligned}
&((\partial_{i-1} \oplus \partial'_{i-1}) \circ \bar{f}_i)(x
\oplus y)=(\partial_{i-1} \oplus \partial'_{i-1})(f_i(x) \oplus
(g_i(x)+f'_i(y)))\\
&=\partial_{i-1}(f_i(x)) \oplus (\partial'_{i-1}(g_i(x)) +
\partial'_{i-1}(f'_i(y)))\\
&=f_{i-1}(\partial_{i-1}(x)) \oplus (g_{i-1}(\partial_{i-1}(x)) +
f'_{i-1}(\partial'_{i-1}(y))\\
&=\bar{f}_{i-1}(\partial_{i-1}(x) \oplus \partial'_{i-1}(y))
=(\bar{f}_{i-1} \circ (\partial_{i-1} \oplus \partial'_{i-1}))(x
\oplus y).
\end{aligned}$$ Hence, $\bar{f}:C \oplus C' \rightarrow C \oplus C'$
is a chain map.

Step 2. Show that $\bar{f}:C \oplus C' \rightarrow C \oplus C'$ is
a quasi-isomorphism.

Notice that the induced homomorphism $\bar{f}_{i*}:H_i(C) \oplus
H_i(C') \rightarrow H_i(C) \oplus H_i(C')$ is defined by
$$\bar{f}_{i*}([x] \oplus [y])=[f_i(x)] \oplus [g_i(x)+f'_i(y)]$$
for all $x \in Z_i(C)$ and $y \in Z_i(C')$ for each $i \in
\{0,\dots,m\}$, where $[x]=x+B_i(C)$, $[f_i(x)]=f_i(x)+B_i(C)$,
$[y]=y+B_i(C')$, and $[g_i(x)+f'_i(y)]=g_i(x)+f'_i(y)+B_i(C')$.

Let $i \in \{0,\dots,m\}$. We claim that $\bar{f}_{i*}$ is an
isomorphism.

Suppose that $x_1,x_2 \in Z_i(C)$ and $y_1,y_2 \in Z_i(C')$ and
$[f_i(x_1)] \oplus [g_i(x_1)+f'_i(y_1)]=[f_i(x_2)] \oplus
[g_i(x_2)+f'_i(y_2)]$. Then $[f_i(x_1)]=[f_i(x_2)]$ and
$[g_i(x_1)+f'_i(y_1)]=[g_i(x_2)+f'_i(y_2)]$. Since $f_{i*}$ is
1-to-1, $[x_1]=[x_2]$, so $[g_i(x_1)]=[g_i(x_2)]$. Hence,
$[f'_i(y_1)]=[f'_i(y_2)]$. Since $f'_{i*}$ is 1-to-1,
$[y_1]=[y_2]$. Therefore, $[x_1] \oplus [y_1]=[x_2] \oplus [y_2]$.

To show that $\bar{f}_{i*}$ is onto, let $x \in Z_i(C)$ and $y \in
Z_i(C')$. Since $f_{i*}$ is onto, there is $a \in Z_i(C)$ such
that $[x]=[f_i(a)]$. Notice that $y - g_i(a) \in Z_i(C')$. Since
$f'_{i*}$ is onto, we can take $b \in Z_i(C')$ such that $[y -
g_i(a)]=[f'_i(b)]$. Hence, $[y]=[g_i(a) + f'_i(b)]$. That is, $[x]
\oplus [y]=[f_i(a)] \oplus [g_i(a) + f'_i(b)]=\bar{f}_{i*}([a]
\oplus [b])$. Hence, $\bar{f}:C \oplus C' \rightarrow C \oplus C'$
is a quasi-isomorphism.

Step 3. Show that $\tau(\bar{f})=\tau(f)\tau(f')$.

Since $\bar{f}$ is a quasi-isomorphism from $C \oplus C'$ to
itself, by Theorem 3.20,
$$\tau(\bar{f})=\prod_{i=0}^m [h_i \oplus h'_i /\bar{f}_{i*}(h_i
\oplus h'_i)]^{(-1)^{i+1}}=\prod_{i=0}^m [h_i \oplus h'_i
/f_{i*}(h_i) \oplus (g_{i*}(h_i) + f'_{i*}(h'_i))]^{(-1)^{i+1}}.$$
Note that for each $i \in \{0,\dots,m\}$, $$(h_i \oplus h'_i
/f_{i*}(h_i) \oplus (g_{i*}(h_i) + f'_{i*}(h'_i)))=(f_{i*}(h_i)
\oplus (g_{i*}(h_i) + f'_{i*}(h'_i))/h_i \oplus h'_i)^{-1}$$ and
$$(f_{i*}(h_i) \oplus (g_{i*}(h_i) + f'_{i*}(h'_i))/h_i \oplus
h'_i)=\begin{pmatrix} (f_{i*}(h_i)/h_i) & (g_{i*}(h_i)/h'_i) \\ 0
& (f'_{i*}(h'_i)/h'_i) \end{pmatrix}.$$ Hence, $$[f_{i*}(h_i)
\oplus (g_{i*}(h_i) + f'_{i*}(h'_i))/h_i \oplus
h'_i]=[f_{i*}(h_i)/h_i][f'_{i*}(h'_i)/h'_i].$$ That is,
$$[h_i \oplus h'_i/f_{i*}(h_i) \oplus (g_{i*}(h_i) +
f'_{i*}(h'_i))]=[h_i/f_{i*}(h_i)][h'_i/f'_{i*}(h'_i)].$$
Therefore, $$\tau(\bar{f})=\prod_{i=0}^m
([h_i/f_{i*}(h_i)][h'_i/f'_{i*}(h'_i)])^{(-1)^{i+1}}=\tau(f)\tau(f').$$

Obviously, if $g=0$, then $\bar{f}=f \oplus f'$. By Lemma 3.6,
$\tau(\bar{f})=\pm\,\tau(f)\tau(f')$. Since $f$ is a
quasi-isomorphism from $C$ to itself and $f'$ is a
quasi-isomorphism from $C'$ to itself, we have
$\tau(\bar{f})=\tau(f)\tau(f')$.

This proves the theorem.
\end{proof}

Now, let us consider the quotient of torsion of
quasi-isomorphisms.

\begin{thm} Let $C$ and $C'$ be based chain complexes of length $m$
over a field $F$, and let $f:C \rightarrow C'$ and $g:C
\rightarrow C'$ be quasi-isomorphisms. Then
$$\frac{\tau(f)}{\tau(g)}=\prod_{i=0}^m
[g_{i*}(h_i)/f_{i*}(h_i)]^{(-1)^{i+1}},$$ where $h_i$ is a basis
for $H_i(C)$ for each $i \in \{0,\dots,m\}$.
\end{thm}

\begin{proof} Suppose that $b_i$, $b_{i-1}$, $b'_i$,
$b'_{i-1}$, $c_i$, and $c'_i$ are bases for $B_i(C)$,
$B_{i-1}(C)$, $B_i(C')$, $B_{i-1}(C')$, $C_i$, and $C'_i$,
respectively, for each $i \in \{0,\dots,m\}$. Then
$$\begin{aligned}
\frac{\tau(f)}{\tau(g)} &=\prod_{i=0}^m
\left(\frac{[(b_ih_i)b_{i-1}/c_i]/[(b'_if_{i*}(h_i))b'_{i-1}/c'_i]}
{[(b_ih_i)b_{i-1}/c_i]/[(b'_ig_{i*}(h_i))b'_{i-1}/c'_i]}
\right)^{(-1)^{i+1}}\\
&=\prod_{i=0}^m \left(\frac{[(b'_ig_{i*}(h_i))b'_{i-1}/c'_i]}
{[(b'_if_{i*}(h_i))b'_{i-1}/c'_i]} \right)^{(-1)^{i+1}}\\
&=\prod_{i=0}^m
[(b'_ig_{i*}(h_i))b'_{i-1}/(b'_if_{i*}(h_i))b'_{i-1}]
^{(-1)^{i+1}}=\prod_{i=0}^m
[g_{i*}(h_i)/f_{i*}(h_i)]^{(-1)^{i+1}}.
\end{aligned}$$
\end{proof}

Note that $\prod_{i=0}^m
[g_{i*}(h_i)/f_{i*}(h_i)]^{(-1)^{i+1}}=\prod_{i=0}^m
[f_{i*}(h_i)/g_{i*}(h_i)]^{(-1)^i}$.

\section{The torsion of dual map of a quasi-isomorphism}

In this section, we introduce the duality theorem for torsion of a
quasi-isomorphism. To prove this theorem, we need the following
statements.

\begin{pro} Let $V$ and $W$ be finite dimensional vector spaces over
a field $F$, and let $f:V \rightarrow W$ be a linear
transformation. Then if $v$ and $w$ be ordered bases for $V$ and
$W$, respectively, then $(f(v)/w)=(f^*(w^*)/v^*)^t$, where
$f^*:W^* \rightarrow V^*$ is the dual map of $f$ and $f(v)$ and
$f^*(w^*)$ are the ordered images of $v$ and $w^*$ under $f$ and
$f^*$, respectively.
\end{pro}

\begin{proof} Suppose that ${\rm dim}_F\,V=m$ and ${\rm
dim}_F\,W=n$ and $v=(v_1,\dots,v_m)$ and $w=(w_1,\dots,w_n)$.
Assume that the matrix $(f(v)/w)$ of $f$ and the matrix
$(f^*(w^*)/v^*)$ of $f^*$ are as follows.
$$(f(v)/w)=
\begin{pmatrix}
a_{11} & a_{12} & \cdots & a_{1n} \\
a_{21} & a_{22} & \cdots & a_{2n} \\
\vdots & \vdots & \cdots & \vdots \\
a_{m1} & a_{m2} & \cdots & a_{mn} \\
\end{pmatrix}
\,\,{\rm and}\,\,
(f^*(w^*)/v^*)=
\begin{pmatrix}
b_{11} & b_{12} & \cdots & b_{1m} \\
b_{21} & b_{22} & \cdots & b_{2m} \\
\vdots & \vdots & \cdots & \vdots \\
b_{n1} & b_{n2} & \cdots & b_{nm} \\
\end{pmatrix}.$$
Let $i \in \{1,\dots,m\}$ and $j \in \{1,\dots,n\}$. Then
$$f(v_i)=\sum_{k=1}^{n}a_{ik}w_k \,\,\,\,\,{\rm and}\,\,\,\,\,
f^*(w_j^*)=\sum_{k=1}^{m}b_{jk}v_k^*.$$ Hence, we have
$$(f^*(w_j^*))(v_i)=(w_j^* \circ f)(v_i)=w_j^*
\left(\sum_{k=1}^{n}a_{ik}w_k
\right)=\sum_{k=1}^{n}a_{ik}w_j^*(w_k)=a_{ij}$$ and
$$\left(\sum_{k=1}^{m}b_{jk}v_k^* \right)(v_i)=
\sum_{k=1}^{m}b_{jk}v_k^*(v_i)=b_{ji}.$$
Therefore, $(f(v)/w)=(f^*(w^*)/v^*)^t$. That is, the matrix
representations are the transpose of each other.
\end{proof}

\begin{pro} Let $V_1$, $V_2$, and $V_3$ be finite dimensional vector
spaces over a field $F$. Then if $f:V_1 \rightarrow V_2$ and
$g:V_2 \rightarrow V_3$ are linear transformations, then $(g \circ
f)^*=f^* \circ g^*$. If $f=0$, then $f^*=0$.
\end{pro}

\begin{proof} Let $h \in V_3^*$. Then $h:V_3 \rightarrow F$ is a
linear functional and $(g \circ f)^*(h) = h \circ (g \circ f)$ and
$(f^* \circ g^*)(h)=f^*(g^*(h))=f^*(h \circ g)=(h \circ g) \circ
f$. Hence, $(g \circ f)^*=f^* \circ g^*$. Also, if $f=0$ and $h
\in V_2^*$, then $f^*(h)=h \circ f=0$, so $f^*=0$.
\end{proof}

Notice that, since the matrix representations $(f(v)/w)$ and
$(f^*(w^*)/v^*)$ are the transpose of each other, they have the
same rank. If $V$ and $W$ have the same dimension, the matrix
representations have the same determinant.

To get a little simpler expression, from time to time, we use
$f^*b^*$ for $f^*(b^*)$ as follows.

\begin{lm} Let $V$ and $W$ be finite dimensional vector spaces over
a field $F$, and let $f:V \rightarrow W$ be a linear
transformation, and let $b$ be a basis for ${\rm Im}\,f$. Then if
$\widetilde{b}$ is a lifting of $b$ by $f$, then
$f^*b^*=\widetilde{b}^*$ and $f^*b^*$ is a basis for ${\rm
Im}\,f^*$.
\end{lm}

\begin{proof} Suppose that ${\rm dim}_F\,V=n$ and
${\rm dim}_F\,{\rm Im}\,f=r$ and $k$ is an ordered basis for ${\rm
Ker}\,f$. Let $k=(k_1,\dots,k_{n-r})$ and $b=(b_1,\dots,b_r)$.
Then $\widetilde{b}=(\widetilde{b}_1,\dots,\widetilde{b}_r)$ and
$(k,\widetilde{b})$ is an ordered basis for $V$. Let $i \in
\{1,\dots,r\}$. Then for each $t \in \{1,\dots,n-r\}$,
$(f^*b_i^*)(k_t)=(b_i^* \circ
f)(k_t)=b_i^*(0)=0=\widetilde{b}_i^*(k_t)$ and for each $j \in
\{1,\dots,r\}$, $(f^*b_i^*)(\widetilde{b}_j)=(b_i^* \circ
f)(\widetilde{b}_j)=b_i^*(b_j)=\delta_{ij}=\widetilde{b}_i^*
(\widetilde{b}_j)$. Hence, $f^*b^*=\widetilde{b}^*$. Since
$f^*b^*$ is a linearly independent subset of ${\rm Im}\,f^*$
containing $r$ elements and ${\rm dim}_F\,{\rm Im}\,f^*={\rm
dim}_F\,{\rm Im}\,f=r$, $f^*b^*$ is a basis for ${\rm Im}\,f^*$.
\end{proof}

\begin{lm} Let $C$ be a chain complex over a field $F$, and let
$\alpha \in C^i$. Then

(1) $\alpha \in Z^i(C)$ if and only if $\alpha|_{B_i(C)}=0$; (2)
$\alpha \in B^i(C)$ if and only if $\alpha|_{Z_i(C)}=0$, where
$Z^i(C)={\rm Ker}\,\{\delta^i:C^i \rightarrow C^{i+1}\}$ and
$B^i(C)={\rm Im}\,\{\delta^{i-1}:C^{i-1} \rightarrow C^i\}$. Also,
$C^i=C_i^*$ and $\delta^i=\partial_i^*$ for each $i$.
\end{lm}

\begin{proof} (1) Let $\alpha \in Z^i(C)$. Then $\delta^i\alpha=0$,
so $\delta^i\alpha=\alpha \circ \partial_i=0$. Hence, $(\alpha
\circ \partial_i)(C_{i+1})=\alpha(B_i(C))=0$. That is,
$\alpha|_{B_i(C)}=0$. Conversely, if $\alpha|_{B_i(C)}=0$, then
$\alpha(B_i(C))=(\alpha \circ \partial_i)(C_{i+1})=0$. Hence,
$\delta^i\alpha=\alpha \circ \partial_i=0$. That is, $\alpha \in
Z^i(C)$. (2) Let $\alpha \in B^i(C)$. Then
$\alpha=\delta^{i-1}\beta$ for some $\beta \in C^{i-1}$. Hence,
$\alpha=\beta \circ \partial_{i-1}$ and $\alpha(z)=(\beta \circ
\partial_{i-1})(z)=\beta(0)=0$ for all $z \in Z_i(C)$,
so $\alpha|_{Z_i(C)}=0$. Conversely, if $\alpha|_{Z_i(C)}=0$, then
${\rm Ker}\,\partial_{i-1}=Z_i(C) \subseteq {\rm Ker}\,\alpha$.
Hence, there is a unique homomorphism $\beta:B_{i-1}(C)
\rightarrow F$ such that $\alpha=\beta \circ \partial_{i-1}$.
Extend $\beta$ to $\widetilde{\beta}:C_{i-1} \rightarrow F$.
Therefore, $\alpha=\widetilde{\beta} \circ
\partial_{i-1}=\delta^{i-1}\widetilde{\beta} \in B^i(C)$.
\end{proof}

\begin{lm} The Universal Coefficient Theorem for a Field.
If $C$ is a chain complex of length $m$ over a field $F$, then
$H^i(C) \cong H_i(C)^*$ for each $i \in \{0,\dots,m\}$, where
$H^i(C)=H_{m-i}(C^*)$ and $H_i(C)^*={\rm Hom}_F(H_i(C),F)$.
\end{lm}

\begin{proof} Let $i \in \{0,\dots,m\}$. Define a function
$\Phi_i:Z^i(C) \rightarrow {\rm Hom}_F(H_i(C),F)$ by
$\Phi_i(\alpha)=(\alpha|_{Z_i(C)})_*$ for each $\alpha \in
Z^i(C)$. We claim that $\Phi_i$ is an epimorphism. Let $\alpha \in
Z^i(C)$. Then $\alpha|_{B_i(C)}=0$, so there is a unique
homomorphism $(\alpha|_{Z_i(C)})_*:H_i(C) \rightarrow F$ such that
$\alpha|_{Z_i(C)}=(\alpha|_{Z_i(C)})_* \circ \pi$, where
$\pi:Z_i(C) \rightarrow H_i(C)$ is the canonical map. Hence,
$\Phi_i$ is well-defined. Let $\alpha, \alpha' \in Z^i(C)$. Then
$\alpha|_{B_i(C)}=\alpha'|_{B_i(C)}=0$. For each $z \in Z_i(C)$,
$((\alpha + \alpha')|_{Z_i(C)})_*([z])=(\alpha +
\alpha')(z)=\alpha(z) +
\alpha'(z)=(\alpha|_{Z_i(C)})_*([z])+(\alpha'|_{Z_i(C)})_*([z])=
((\alpha|_{Z_i(C)})_* + (\alpha'|_{Z_i(C)})_*)([z])$. Also, for
each $r \in F$,
$((r\alpha)|_{Z_i(C)})_*([z])=(r\alpha)(z)=r\alpha(z)=
r(\alpha|_{Z_i(C)})_*([z])=(r(\alpha|_{Z_i(C)})_*)([z])$, hence,
$\Phi_i$ is a homomorphism. To show that $\Phi_i$ is onto, let
$\beta \in {\rm Hom}_F(H_i(C),F)$. Then $\beta \circ \pi:Z_i(C)
\rightarrow F$ is a homomorphism and $(\beta \circ
\pi)|_{B_i(C)}=0$. Extend $\beta \circ \pi$ to $\widetilde{\beta
\circ \pi}:C_i \rightarrow F$. Since $\beta \circ \pi=(\beta \circ
\pi)_* \circ \pi$ and $\pi$ is onto, $\beta=(\beta \circ \pi)_*$.
Hence, $\beta=\Phi_i(\widetilde{\beta \circ \pi})$ and
$\widetilde{\beta \circ \pi} \in Z^i(C)$. That is, $\Phi_i$ is
onto. Next, we show that $\Phi_i|_{B^i(C)}=0$. Suppose that
$\alpha \in B^i(C)$. Then $\alpha|_{Z_i(C)}=0$ by Lemma 4.4, so
$\Phi_i(\alpha)=(\alpha|_{Z_i(C)})_*=0_*=0$. Therefore, we have
the induced homomorphism $\Phi_{i*}:H^i(C) \rightarrow {\rm
Hom}_F(H_i(C),F)$ which is onto. By the universal coefficient
theorem, we have the short exact sequence
$$\begin{CD}
0 @>>> {\rm Ext}(H_{i-1}(C),F) @>>> H^i(C) @>>> {\rm
Hom}_F(H_i(C),F) @>>> 0
\end{CD}\,.$$
Since $F$ is a field, ${\rm Ext}(H_{i-1}(C),F)=0$. Therefore,
$H^i(C) \cong H_i(C)^*$.
\end{proof}

\begin{lm} Let $C$ and $C'$ be chain complexes of length $m$ over
a field $F$ such that $C_i$ and $C'_i$ are finite dimensional
vector spaces for each $i \in \{0,\dots,m\}$, and let $f:C
\rightarrow C'$ be a quasi-isomorphism. Then for each $i \in
\{0,\dots,m\}$,

(1) if $h_i$ is a basis for $H_i(C)$, then $[\widetilde{h_i}^*]$
is a basis for $H^i(C)$;

(2) $(f_i^*)_*([f_i(\widetilde{h_i})^*])=[\widetilde{h_i}^*]$.
\end{lm}

\begin{proof} Suppose that ${\rm dim}\,B_i(C)=r$ and
${\rm dim}\,H_i(C)=s$ and ${\rm dim}\,B_{i-1}(C)=t$ and
$b_i=(b_{i1},\dots,b_{ir})$ and $h_i=(h_{i1},\dots,h_{is})$ and
$b_{i-1}=(b_{i-11},\dots,b_{i-1t})$. Notice that
$(b_i,\widetilde{h_i},\widetilde{b_{i-1}})$ is an ordered basis
for $C_i$. Hence,
$(b_i^*,\widetilde{h_i}^*,\widetilde{b_{i-1}}^*)$ is an ordered
basis for $C_i^*$. In particular, $\widetilde{b_{i-1}}^*$ is an
ordered basis for $B^i(C)$.

(1) For each $j \in \{1,\dots,s\}$ and $k \in \{1,\dots,r\}$,
$\widetilde{h_i}_j^*(b_{ik})=0$, that is,
$\widetilde{h_i}_j^*|_{B_i(C)}=0$. Hence, by Lemma 4.4,
$\widetilde{h_i}^* \subset Z^i(C)$. Suppose that
$r_1[\widetilde{h_i}_1^*]+\cdots+r_s[\widetilde{h_i}_s^*]=B^i(C)$.
Then $r_1\widetilde{h_i}_1^*+\cdots+r_s\widetilde{h_i}_s^* \in
B^i(C)$. Hence, by Lemma 4.4,
$(r_1\widetilde{h_i}_1^*+\cdots+r_s\widetilde{h_i}_s^*)|_{Z_i(C)}=0$.
In particular,
$(r_1\widetilde{h_i}_1^*+\cdots+r_s\widetilde{h_i}_s^*)
(\widetilde{h_i}_j)=r_j=0$ for each $j \in \{1,\dots,s\}$. Since
${\rm dim}\,H^i(C)={\rm dim}\,H_i(C)=s$, $[\widetilde{h_i}^*]$ is
a basis for $H^i(C)$.

(2) Note that $f_{i*}(h_i)=[f_i(\widetilde{h_i})]$, where
$f_{i*}(h_i)=(f_{i*}(h_{i1}),\dots,f_{i*}(h_{is}))$ and
$[f_i(\widetilde{h_i})]=([f_i(\widetilde{h_i}_1)],\dots,
[f_i(\widetilde{h_i}_s)])$. Since
$(f_i^*)_*([f_i(\widetilde{h_i})^*])=[f_i^*(f_i(\widetilde{h_i})^*)]$,
it suffices to show that
$f_i^*(f_i(\widetilde{h_i}_j)^*)-\widetilde{h_i}_j^* \in B^i(C)$
for each $j \in \{1,\dots,s\}$. We claim that
$(f_i(\widetilde{h_i}_j)^* \circ
f_i)|_{Z_i(C)}=\widetilde{h_i}_j^*|_{Z_i(C)}$ for each $j \in
\{1,\dots,s\}$. Remark that $(b_i,\widetilde{h_i})$ is a basis for
$Z_i(C)$. Let $j \in \{1,\dots,s\}$ and $k \in \{1,\dots,r\}$ and
$l \in \{1,\dots,s\}$. Then $\widetilde{h_i}_j^*(b_{ik})=0$. since
$f_{i*}$ is an isomorphism, $f_i(\widetilde{h_i}_j) \notin
B_i(C')$, but $f_i(b_{ik}) \in B_i(C')$. Hence,
$f_i(\widetilde{h_i}_j)^*(f_i(b_{ik}))=0$. Also,
$\widetilde{h_i}_j^*(\widetilde{h_i}_l)=\delta_{jl}$ and
$f_i(\widetilde{h_i}_j)^*(f_i(\widetilde{h_i}_l))=\delta_{jl}$
since $f_{i*}$ is an isomorphism. Hence, $(f_i(\widetilde{h_i})^*
\circ f_i)|_{Z_i(C)}=\widetilde{h_i}^*|_{Z_i(C)}$. Therefore,
$(f_i^*)_*([f_i(\widetilde{h_i})^*])=[\widetilde{h_i}^*]$.
\end{proof}

Now we introduce the torsion of dual map of a quasi-isomorphism.
It turns out that the torsion of dual map of a quasi-isomorphism
depends on the length of chain complexes which is domain or range
of the quasi-isomorphism.

\begin{thm} Let $C$ and $C'$ be based chain complexes of length $m$
over a field $F$ such that $C_i$ and $C'_i$ are finite dimensional
vector spaces for each $i \in \{0,\dots,m\}$, and let $f:C
\rightarrow C'$ be a quasi-isomorphism. Then $f^*:C'^* \rightarrow
C^*$ is a quasi-isomorphism and $\tau(f^*)=\pm\,\tau(f)^{(-1)^m}$.
\end{thm}

\begin{proof} For each $i \in \{0,\dots,m-1\}$, $\partial_i^* \circ
f_i^*=(f_i \circ \partial_i)^*=(\partial'_i \circ
f_{i+1})^*=f_{i+1}^* \circ {\partial'}_i^*$. Hence, $f^*:C'^*
\rightarrow C^*$ is a chain map. Since $f_{i*}:H_i(C) \rightarrow
H_i(C')$ is an isomorphism, $(f_{i*})^*:H_i(C')^* \rightarrow
H_i(C)^*$ is an isomorphism. In the proof of the universal
coefficient theorem for a field (Lemma 4.5), we defined the
isomorphism $\Phi_{i*}:H^i(C) \rightarrow H_i(C)^*$. Define the
isomorphism $\theta_{i*}:H^i(C') \rightarrow H_i(C')^*$ by the
same argument as for $C$. Now we claim that
$(f_i^*)_*=\Phi_{i*}^{-1} \circ (f_{i*})^* \circ \theta_{i*}$. Let
$\alpha \in Z^i(C')$. Then $[\alpha] \in H^i(C')$ and
$\theta_{i*}([\alpha])=(\alpha|_{Z_i(C')})_*$ and
$(f_{i*})^*(\theta_{i*}([\alpha]))=(\alpha|_{Z_i(C')})_* \circ
f_{i*}=(\alpha|_{Z_i(C')} \circ f_i)_*=((\alpha \circ
f_i)|_{Z_i(C)})_*$. Hence, $(\Phi_{i*}^{-1} \circ (f_{i*})^* \circ
\theta_{i*})([\alpha])=[\alpha \circ
f_i]=[f_i^*(\alpha)]=(f_i^*)_*([\alpha])$. Therefore,
$(f_i^*)_*:H^i(C') \rightarrow H^i(C)$ is an isomorphism, hence,
$f^*:C'^* \rightarrow C^*$ is a quasi-isomorphism.

Now we consider the torsion $\tau(f^*)$ of the dual
quasi-isomorphism $f^*$ of $f$.

Note that for each $i \in \{0,\dots,m\}$, $(f^*)_i=(f_{m-i})^*$,
$(C^*)_i=(C_{m-i})^*$, $(C'^*)_i=(C'_{m-i})^*$,
$(\partial^*)_i=(\partial_{m-i-1})^*$, and
$(\partial'^*)_i=(\partial'_{m-i-1})^*$.

Let us use the convention as follows. For each $i \in
\{0,\dots,m\}$, $(f_i)^*=f_i^*$, $(C_i)^*=C_i^*$,
$(C'_i)^*={C'}_i^*$, $(\partial_i)^*=\partial_i^*$, and
$(\partial'_i)^*={\partial'}_i^*$. Then

$$\begin{aligned}
\tau(f^*) &=\prod_{j=0}^m
\left(\frac{[((b'^*)_j(h'^*)_j)(b'^*)_{j-1}/(c'^*)_j]}
{[((b^*)_j(f^*)_{j*}((h'^*)_j))(b^*)_{j-1}/(c^*)_j]}
\right)^{(-1)^{j+1}}\\
&=\prod_{i=0}^m
\left(\frac{[((b'^*)_{m-i}(h'^*)_{m-i})(b'^*)_{m-i-1}/(c'^*)_{m-i}]}
{[((b^*)_{m-i}(f^*)_{m-i*}((h'^*)_{m-i}))(b^*)_{m-i-1}/(c^*)_{m-i}]}
\right)^{(-1)^{m-i+1}}\\
&=\prod_{i=0}^m
\left(\frac{[({\partial'}_{i-1}^*{b'}_{i-1}^*,\,f_i(\widetilde{h_i})^*,\,
\widetilde{{\partial'}_i^*{b'}_i^*})/{c'}_i^*]}
{[(\partial_{i-1}^*b_{i-1}^*,\,f_i^*f_i(\widetilde{h_i})^*,\,
\widetilde{\partial_i^*b_i^*})/c_i^*]}
\right)^{(-1)^{m-i+1}}\,\,\,\,\,\textrm{By}\,\textrm{Lemma}\,3.4,\,4.6.\\
&=\prod_{i=0}^m
\left(\frac{[(\widetilde{b'_{i-1}}^*,\,f_i(\widetilde{h_i})^*,\,
{b'}_i^*)/{c'}_i^*]}
{[(\widetilde{b_{i-1}}^*,\,f_i^*f_i(\widetilde{h_i})^*,\,b_i^*)/c_i^*]}
\right)^{(-1)^{m-i+1}}\,\,\,\,\,\,\,\,\,\,\,\,\,\,\,\,\,\,\,\,\,
\textrm{By}\,\textrm{Lemma}\,4.3.\\
&=\prod_{i=0}^m
\left(\frac{[(\widetilde{b'_{i-1}}^*,\,f_i(\widetilde{h_i})^*,\,
{b'}_i^*)/{c'}_i^*]}
{[(\widetilde{b_{i-1}}^*,\,\widetilde{h_i}^*,\,b_i^*)/c_i^*]}
\right)^{(-1)^{m-i+1}}\,\,\,\,\,\,\,\,\,\,\,\,\,\,\,\,\,\,\,\,\,\,\,\,
\textrm{By}\,\textrm{Lemma}\,4.6.\\
&=\pm\,\prod_{i=0}^m
\left(\frac{[({b'}_i^*,\,f_i(\widetilde{h_i})^*,\,\widetilde{b'_{i-1}}^*)
/{c'}_i^*]}
{[(b_i^*,\,\widetilde{h_i}^*,\,\widetilde{b_{i-1}}^*)/c_i^*]}
\right)^{(-1)^{m-i+1}}\\
&=\pm\,\prod_{i=0}^m
\left(\frac{[(b'_i,\,f_i(\widetilde{h_i}),\,\widetilde{b'_{i-1}})
/c'_i]^{-1}}
{[(b_i,\,\widetilde{h_i},\,\widetilde{b_{i-1}})/c_i]^{-1}}
\right)^{(-1)^{m-i+1}}\,\,\,\,\,\,\,\,\,\,\,\,\,\,\,\,\,\,\,\,
\textrm{By}\,\textrm{Proposition}\,4.1\\
&=\pm\,\prod_{i=0}^m
\left(\frac{[(b_i,\,\widetilde{h_i},\,\widetilde{b_{i-1}})/c_i]}
{[(b'_i,\,f_i(\widetilde{h_i}),\,\widetilde{b'_{i-1}})/c'_i]}
\right)^{(-1)^{i+1}(-1)^m}\\
&=\pm\,\left(\prod_{i=0}^m \left(\frac{[(b_ih_i)b_{i-1}/c_i]}
{[(b'_if_{i*}(h_i))b'_{i-1}/c'_i]}
\right)^{(-1)^{i+1}}\right)^{(-1)^m}=\pm\,\tau(f)^{(-1)^m}.
\end{aligned}$$
\end{proof}

\bigskip
\centerline{\epsfxsize=5 in \epsfbox{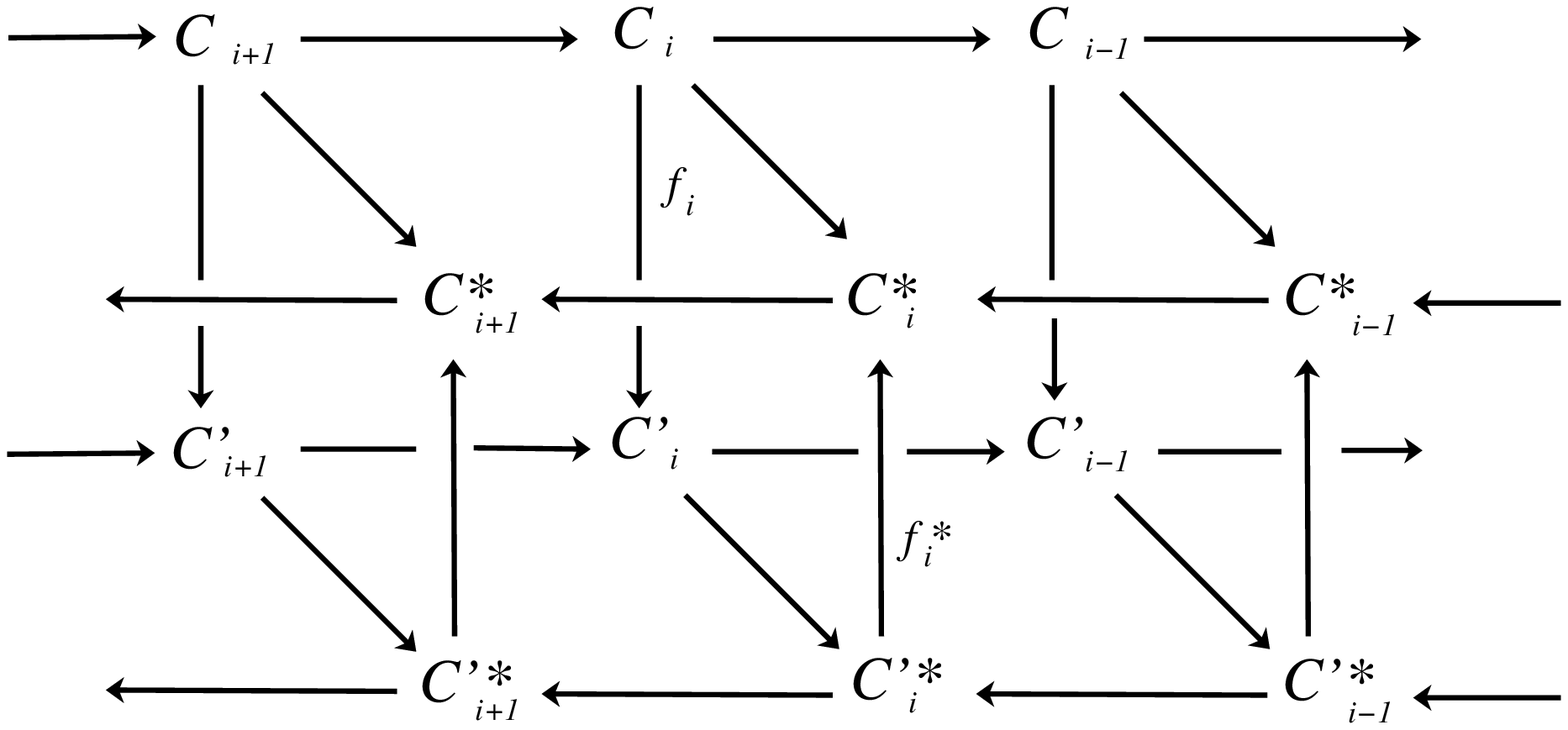}}
\medskip
\centerline{\small Figure 1. The diagram of a quasi-isomorphism
and its dual map.}
\bigskip

As the same idea that we determine the sign in Lemma 3.6, we have
the exact sign for the torsion of dual map.

For each $i \in \{0,\dots,m\}$, let
$$x_i={\rm dim}_F\,B_i(C), \,\,\,x'_i={\rm dim}_F\,B_i(C'),
\,\,\,y_i={\rm dim}_F\,H_i(C).$$ Then
$$\tau(f^*)=\frac{(-1)^{\sum_{i=0}^m(x'_i(x'_{i-1}+y_i) +
x'_{i-1}y_i)}}{(-1)^{\sum_{i=0}^m(x_i(x_{i-1}+y_i) +
x_{i-1}y_i)}}\,\tau(f)^{(-1)^m}.$$ Therefore,
$$\tau(f^*)=(-1)^{\sum_{i=0}^m[(x'_i(x'_{i-1}+y_i) +
x'_{i-1}y_i)-(x_i(x_{i-1}+y_i) +
x_{i-1}y_i)]}\,\tau(f)^{(-1)^m}.$$

\section{torsion of quasi-isomorphisms between free chain complexes}

In this section, we generalize our torsion so that we define the
torsion of a quasi-isomorphism between free chain complexes.

Suppose that $R$ is an associative ring with $1 \ne 0$ which has
invariant dimension property, that is, $m=n$ if and only if $R^m
\cong R^n,$ where $m$ and $n$ are nonnegative integers and $R^m$
and $R^n$ are direct sums of $R$. For each $n \in \mathbb N$, let
${\rm GL}(n,R)$ be the group of $n \times n$ invertible matrices
over $R$, called the $n$-general linear group over $R$. We can
identify each $A \in {\rm GL}(n,R)$ with the matrix
$$\begin{pmatrix} A & 0 \\ 0 & 1 \end{pmatrix} \in {\rm
GL}(n+1,R)$$ so that we consider $${\rm GL}(1,R) \subset {\rm
GL}(2,R) \subset \cdots.$$ ${\rm GL}(R)=\bigcup_{n \in \mathbb
N}{\rm GL}(n,R)$ is called the infinite general linear group over
$R$.

\begin{nota} $$K_1(R)={\rm GL}(R)/[{\rm GL}(R),{\rm GL}(R)],$$
where $[{\rm GL}(R),{\rm GL}(R)]$ is the commutator subgroup of
${\rm GL}(R)$.
\end{nota}

For the time being, let us assume that a ring $R$ has $1 \ne 0$
and invariant dimension property, but need not be commutative.
Consider a based chain complex $C$ over $R$ such that $C_i$ is a
based free $R$-module of finite rank for each $i$. We call such a
chain complex a based free chain complex over $R$.

Suppose that $C$ and $C'$ are based free chain complexes of length
$m$ over $R$ and $f:C \rightarrow C'$ is a quasi-isomorphism. If
$B_i(C)$, $B_i(C')$, and $H_i(C)$ are free $R$-modules for all $i
\in \{0,\dots,m\}$, then we define the torsion $\tau(f)$ of $f:C
\rightarrow C'$ by
$$\tau(f)=\left[\prod_{i=0}^m
\left(\frac{((b_ih_i)b_{i-1}/c_i)}{((b'_if_{i*}(h_i))b'_{i-1}/c'_i)}
\right)^{(-1)^{i+1}}\right] \in K_1(R),$$ where $b_i$, $b_{i-1}$,
$b'_i$, $b'_{i-1}$, $c_i$, $c'_i$, and $h_i$ are bases for
$B_i(C)$, $B_{i-1}(C)$, $B_i(C')$, $B_{i-1}(C')$, $C_i$, $C'_i$,
and $H_i(C)$, respectively, for each $i \in \{0,\dots,m\}$. Also,
$\frac{1}{A}$ means the inverse $A^{-1}$ of $A$ and $[A]$ means
the abelianized class $A[{\rm GL}(R),{\rm GL}(R)]$ of $A$ for $A
\in {\rm GL}(R)$.

By the same idea as used in the proof of Lemma 3.4, we can show
that $\tau(f)$ is well-defined. Note that if $R$ is a commutative
ring with $1$, then the determinant ${\rm det}:K_1(R) \rightarrow
R-\{0\}$ is a surjective group homomorphism. In particular, if $R$
is a field, then ${\rm det}:K_1(R) \rightarrow R-\{0\}$ is an
isomorphism. See \cite{T}.

Therefore, when $R$ is a field, we can identify above definition
with Definition 3.3 by this isomorphism.

In general, even though $C$ is a free chain complex, its boundary
and homology modules need not be free. We show that the torsion
$\tau(f)$ of a quasi-isomorphism $f:C \rightarrow C'$ is defined
if each homology module is a summand of a free module.

\begin{lm} Let $F$ be a free $R$-module, and let $m \in \mathbb N$
and $i \in \{0,\dots,m-1\}$. Then if $C$ is the chain complex of
length $m$ over $R$ such that $C_i=C_{i+1}=F$ and $C_k=0$ if $k
\ne i, \,i+1$ and $\partial_i=I_F$ and $\partial_k=0$ if $k \ne
i$, then $C$ is acyclic.
\end{lm}

\begin{proof} We see that ${\rm Im}\,\partial_i=F$ and ${\rm Ker}\,
\partial_{i-1}=F$. If $j \ne i$, then ${\rm Im}\,\partial_j=0$ and
${\rm Ker}\,\partial_{j-1}=0$. Hence, ${\rm Im}\,\partial_k={\rm
Ker}\,\partial_{k-1}$ for each $k \in \{0,\dots,m\}$. That is, $C$
is acyclic.
\end{proof}
$$\begin{CD}
\cdots @>>> 0 @>>> F @>I_{F}>> F @>>> 0 @>>> \cdots
\end{CD}$$\\

\begin{nota} The chain complex in Lemma 5.2 is denoted by $C(F,i)$
for each $i \in \{0,\dots,m-1\}$.
\end{nota}

\begin{lm} If $C$ and $C'$ are acyclic chain complexes of length
$m$ over $R$, then $C \oplus C'$ is an acyclic chain complex of
length $m$ over $R$.
\end{lm}

\begin{proof} Let $i \in \{0,\dots,m\}$. Since $C$ and $C'$ are acyclic,
${\rm Im}\,\partial_i={\rm Ker}\,\partial_{i-1}$ and ${\rm
Im}\,\partial'_i={\rm Ker}\,\partial'_{i-1}$. Hence, ${\rm
Im}\,\partial_i \oplus {\rm Im}\,\partial'_i={\rm
Ker}\,\partial_{i-1} \oplus {\rm Ker}\,\partial'_{i-1}$. Note that
${\rm Im} (\partial_i \oplus \partial'_i)={\rm Im}\,\partial_i
\oplus {\rm Im}\,\partial'_i$ and ${\rm Ker}(\partial_{i-1} \oplus
\partial'_{i-1})={\rm Ker}\,\partial_{i-1} \oplus {\rm
Ker}\,\partial'_{i-1}$. Therefore, $C \oplus C'$ is acyclic.
\end{proof}

By Lemma 5.2 and 5.4, we have the following statement immediately
which plays an important role for our generalization.

\begin{lm} Let $C$ be an acyclic based free chain complex of length
$m \geq 1$ over $R$, and let $F$ be a free based $R$-module. Then
for each $i \in \{0,\dots,m-1\}$, $C \oplus C(F_i,i)$ is an
acyclic based free chain complex of length $m$ over $R$.
Furthermore, $C \oplus \bigoplus_{i=0}^{m-1} C(F_i,i)$ is also an
acyclic based free chain complex of length $m$ over $R$.
\end{lm}

\begin{df} An $R$-module $M$ is said to be stably free if there is
a free $R$-module $F$ of finite rank such that $M \oplus F$ is
free.
\end{df}

Note that zero $R$-module is stably free.

\begin{lm} Let $C$ be a based free chain complex of length $m$
over $R$. Then if $H_i(C)$ is stably free for all $i \in
\{0,\dots,m\}$, then $Z_i(C)$ and $B_i(C)$ are stably free for all
$i \in \{0,\dots,m\}$.
\end{lm}

\begin{proof} We prove this by induction. Suppose that $H_i(C)$ is
stably free for all $i \in \{0,\dots,m\}$, say, $H_i(C) \oplus
F^H_i$ is free for some free $R$-module $F^H_i$. Remark that the
direct sum of 2 short exact sequences
$$\begin{CD}
0 @>>> B_i(C) @>\subseteq>> Z_i(C) @>\pi>> H_i(C) @>>> 0
\end{CD}$$ and
$$\begin{CD}
0 @>>> 0 @>>> F^H_i @>I_{F^H_i}>> F^H_i @>>> 0
\end{CD}$$ is the short exact sequence
$$\begin{CD}
0 @>>> B_i(C) \oplus 0 @>>> Z_i(C) \oplus F^H_i @>>> H_i(C) \oplus
F^H_i @>>> 0
\end{CD}\,.$$
Since $Z_0(C)=C_0$, $Z_0(C)$ is free, so $Z_0(C)$ is stably free.
Since $H_0(C) \oplus F^H_0$ is free, $Z_0(C) \oplus F^H_0=B_0(C)
\oplus 0 \oplus H_0(C) \oplus F^H_0=B_0(C) \oplus H_0(C) \oplus
F^H_0$, hence, $B_0(C)$ is stably free. Suppose that $B_{i-1}(C)$
is stably free for each $i \in \{1,\dots,m\}$, say, $B_{i-1}(C)
\oplus F^B_{i-1}$ is free for some free $R$-module $F^B_{i-1}$. To
show that $B_i(C)$ is stably free, consider 3 short exact
sequences
$$\begin{CD}
0 @>>> Z_i(C) @>\subseteq>> C_i @>\partial_{i-1}>> B_{i-1}(C) @>>>
0
\end{CD}$$ and
$$\begin{CD}
0 @>>> 0 @>>> F^B_{i-1} @>I_{F^B_{i-1}}>> F^B_{i-1} @>>> 0
\end{CD}$$ and
$$\begin{CD}
0 @>>> F^H_i @>I_{F^H_i}>> F^H_i @>>> 0 @>>> 0
\end{CD}\,.$$ Then we have the short exact sequence
$$0 \longrightarrow Z_i(C) \oplus 0 \oplus F^H_i \longrightarrow
C_i \oplus F^B_{i-1} \oplus F^H_i \longrightarrow B_{i-1}(C)
\oplus F^B_{i-1} \oplus 0 \longrightarrow 0$$ which is the direct
sum of those 3 short exact sequences. Hence,
$$\begin{aligned}
C_i \oplus F^B_{i-1} \oplus F^H_i
&=Z_i(C) \oplus F^H_i \oplus
B_{i-1}(C) \oplus F^B_{i-1}\\
&=B_i(C) \oplus H_i(C) \oplus F^H_i \oplus B_{i-1}(C) \oplus
F^B_{i-1}.
\end{aligned}$$
Since $H_i(C) \oplus F^H_i \oplus B_{i-1}(C) \oplus F^B_{i-1}$ is
free, $B_i(C)$ is stably free. Hence, $B_i(C)$ are stably free for
all $i \in \{0,\dots,m\}$.

Let $i \in \{0,\dots,m\}$. Suppose that $F^B_i$ and $B_i(C) \oplus
F^B_i$ are free. Then we have the exact sequence
$$0 \longrightarrow B_i(C) \oplus 0 \oplus F^B_i \longrightarrow
Z_i(C) \oplus F^H_i \oplus F^B_i \longrightarrow H_i(C) \oplus
F^H_i \oplus 0 \longrightarrow 0.$$ Since $H_i(C) \oplus F^H_i$ is
free, $Z_i(C) \oplus F^H_i \oplus F^B_i=B_i(C) \oplus F^B_i \oplus
H_i(C) \oplus F^H_i$ and $Z_i(C) \oplus F^H_i \oplus F^B_i$ is
free, so $Z_i(C)$ is stably free. This proves the lemma.
\end{proof}

Let $C$ be a based free chain complex of length $m$ over $R$.
Suppose that $H_i(C)$ is stably free for all $i \in
\{0,\dots,m\}$. Then for each $i \in \{0,\dots,m\}$, by Lemma 5.7,
there are free $R$-modules $F^B_i$ and $F^H_i$ such that $B_i(C)
\oplus F^B_i$ and $H_i(C) \oplus F^H_i$ are free. Notice that $C
\oplus \bigoplus_{i=0}^{m-1} C(F^B_i,i) \oplus \bigoplus_{i=0}^m
F^H_i$ is a free chain complex of length $m$ over $R$ whose
boundary modules and homology modules are free. In particular,
$$H_i(C \oplus \bigoplus_{i=0}^{m-1} C(F^B_i,i) \oplus
\bigoplus_{i=0}^m F^H_i)=H_i(C) \oplus F^H_i$$ for each $i \in
\{0,\dots,m\}$, where $F^H_i$ is the free chain complex $C'$ of
length $m$ over $R$ such that the $C'_i=F^H_i$ and $C'_j=0$ if $j
\ne i$ and $\partial'_i=0$ for all $i \in \{0,\dots,m-1\}$.

Now, we generalize Definition 3.3 to the torsion of a
quasi-isomorphism between free chain complexes.

\begin{df} Let $C$ and $C'$ be based free chain complexes of length $m$
over $R$, and let $f:C \rightarrow C'$ be a quasi-isomorphism.
Suppose that for each $i \in \{0,\dots,m\}$, $F^B_i$, $F^H_i$,
$F'^B_i$, and $F'^H_i$  are free $R$-modules such that $B_i(C)
\oplus F^B_i$, $H_i(C) \oplus F^H_i$, $B_i(C') \oplus F'^B_i$, and
$H_i(C') \oplus F'^H_i$ are free. Define the torsion $\tau(f)$ by
$$\tau(f)=\pm\,\tau(f \oplus S) \in K_1(R),$$ where $S$ is the identity
chain isomorphism from
$$\bigoplus_{i=0}^{m-1} C(F^B_i,i) \oplus \bigoplus_{i=0}^m F^H_i
\oplus \bigoplus_{i=0}^{m-1} C(F'^B_i,i) \oplus \bigoplus_{i=0}^m
F'^H_i$$ to $$\bigoplus_{i=0}^{m-1} C(F^B_i,i) \oplus
\bigoplus_{i=0}^m F^H_i \oplus \bigoplus_{i=0}^{m-1} C(F'^B_i,i)
\oplus \bigoplus_{i=0}^m F'^H_i.$$ The identity chain isomorphism
$S$ is called a stabilizer of $f$.
\end{df}

Notice that, by definition, $\tau(S)=1$ and $\tau(f \oplus
S)=\pm\,\tau(S \oplus f)$.

\begin{lm} Let $C$ and $C'$ be based free chain complexes of length
$m$ over $R$, and let $f:C \rightarrow C'$ be a quasi-isomorphism.
Then if $H_i(C)$ and $H_i(C')$ are stably free for all $i \in
\{0,\dots,m\}$ and $S_1$ and $S_2$ are stabilizers of $f$, then
$\tau(f \oplus S_1)=\pm\,\tau(f \oplus S_2)$. That is, $\tau(f)$
is independent of the choice of stabilizer of $f$ upto sign.
\end{lm}

\begin{proof} Suppose that $f:C \rightarrow C'$ is a quasi-isomorphism
between based free chain complexes of length $m$ over $R$ and
$H_i(C)$ and $H_i(C')$ are stably free for all $i \in
\{0,\dots,m\}$ and $S_1$ and $S_2$ are stabilizers of $f$. Then
$$\tau((f \oplus S_1) \oplus S_2)=\pm\,\tau(f \oplus
S_1)\tau(S_2)=\pm\,\tau(f \oplus S_1)$$ and
$$\begin{aligned}
&\tau((f \oplus S_1) \oplus S_2)=\pm\,\tau(S_2 \oplus (f \oplus
S_1))=\pm\,\tau((S_2 \oplus f) \oplus S_1))\\
&=\pm\,\tau(S_2 \oplus f)\tau(S_1)=\pm\,\tau(S_2 \oplus
f)=\pm\,\tau(f \oplus S_2).
\end{aligned}$$ Therefore, $\tau(f \oplus S_1)=\pm\,\tau(f \oplus S_2)$.
\end{proof}

Next, we briefly introduce another extension of our torsion of
quasi-isomorphisms which generalizes Turaev's Theorem \cite{T2}
for the torsion of chain complexes whose rank of homology is zero.

\begin{thm} Turaev \cite{T2}. Let $R$ be a Noetherian unique factorization
domain, and let $C$ be a based free chain complex of length $m$
over $R$ such that ${\rm rank}\,H_i(C)=0$ for all $i \in
\{0,\dots,m\}$. Then $\widetilde{C}=\widetilde{R} \otimes_R C$ is
a based acyclic chain complex of vector spaces over
$\widetilde{R}$ and $\tau(C)$ is defined by $\tau(\widetilde{C})$.
Furthermore,
$$\tau(\widetilde{C})=\prod_{i=0}^m \left({\rm
ord}\,H_i(C)\right)^{(-1)^{i+1}},$$ where ${\rm ord}\,H_i(C)$ is
the $0$-th Alexander polynomial $\Delta_0(H_i(C))$ of $H_i(C)$ for
each $i \in \{0,\dots,m\}$.
\end{thm}

Suppose that $f:C \rightarrow C'$ is a chain map which is not a
necessarily quasi-isomorphism, but $\widetilde{f}={\rm id} \otimes
f:\widetilde{R} \otimes_R C \longrightarrow \widetilde{R}
\otimes_R C'$ can be a quasi-isomorphism. In this case, we can
define the torsion of $f$.

\begin{df} Let $R$ be a Noetherian unique factorization
domain, and let $C$ and $C'$ be based free chain complexes of
length $m$ over $R$, and let $f:C \rightarrow C'$ be a chain map
such that $\widetilde{f}={\rm id} \otimes f:\widetilde{R}
\otimes_R C \longrightarrow \widetilde{R} \otimes_R C'$ is a
quasi-isomorphism. Then the torsion $\tau(f)$ of $f$ is defined by
$\tau(f)=\tau(\widetilde{f})$.
\end{df}

This generalized torsion of chain maps has a possible application
to link theory.

$$\begin{CD}
\widetilde{R} \otimes_R C_i @>{\rm id} \otimes \partial_{i-1}>>
\widetilde{R} \otimes_R C_{i-1}\\
@V{\rm id} \otimes f_iVV @VV{\rm id} \otimes f_{i-1}V\\
\widetilde{R} \otimes_R C'_i @>{\rm id} \otimes \partial'_{i-1}>>
\widetilde{R} \otimes_R C'_{i-1}
\end{CD}$$

\begin{cor} Let $R$ be a Noetherian unique factorization
domain, and let $C$ and $C'$ be based free chain complexes of
length $m$ over $R$ such that ${\rm rank}\,H_i(C)=0$ and ${\rm
rank}\,H_i(C')=0$ for all $i \in \{0,\dots,m\}$, and let $f:C
\rightarrow C'$ be a chain map such that $\widetilde{f}={\rm id}
\otimes f:\widetilde{R} \otimes_R C \longrightarrow \widetilde{R}
\otimes_R C'$ is a quasi-isomorphism. Then
$$\tau(f)=\tau(\widetilde{f})=\frac{\tau(\widetilde{C})}
{\tau(\widetilde{C'})}=\prod_{i=0}^m \left(\frac{{\rm
ord}\,H_i(C)}{{\rm ord}\,H_i(C')}\right)^{(-1)^{i+1}},$$ where
$\widetilde{C}=\widetilde{R} \otimes_R C$ and
$\widetilde{C'}=\widetilde{R} \otimes_R C'$.
\end{cor}

\section{Examples of quasi-isomorphisms}

In this section, we introduce a few concrete examples of torsion
of quasi-isomorphism. Notice that in all these examples, the
quasi-isomorphisms are from a chain complex to itself. So by
Theorem 3.20, their torsion can be calculated from their actions
on homology. Nevertheless, we want to show the calculation from
definition so that the reader may get familiar with the
construction.

Note that a vector space $C$ can be regarded as a chain complex $0
\rightarrow C \rightarrow 0$ with length $0$ and a bijective
linear transformation $f:C \rightarrow C$ is a quasi-isomorphism.
In this case, the torsion of $f$ is exactly same as ${\rm det}\,f$
(Theorem 3.20).

$$\begin{CD}
0 @>\partial_1>> C_1 @>\partial_0>> C_0 @>\partial_{-1}>> 0\\
@VVV @Vf_1VV @Vf_0VV @VVV\\
0 @>\partial_1>> C_1 @>\partial_0>> C_0 @>\partial_{-1}>> 0\\
\end{CD}$$\\

\bigskip
\centerline{\epsfxsize=5.8 in \epsfbox{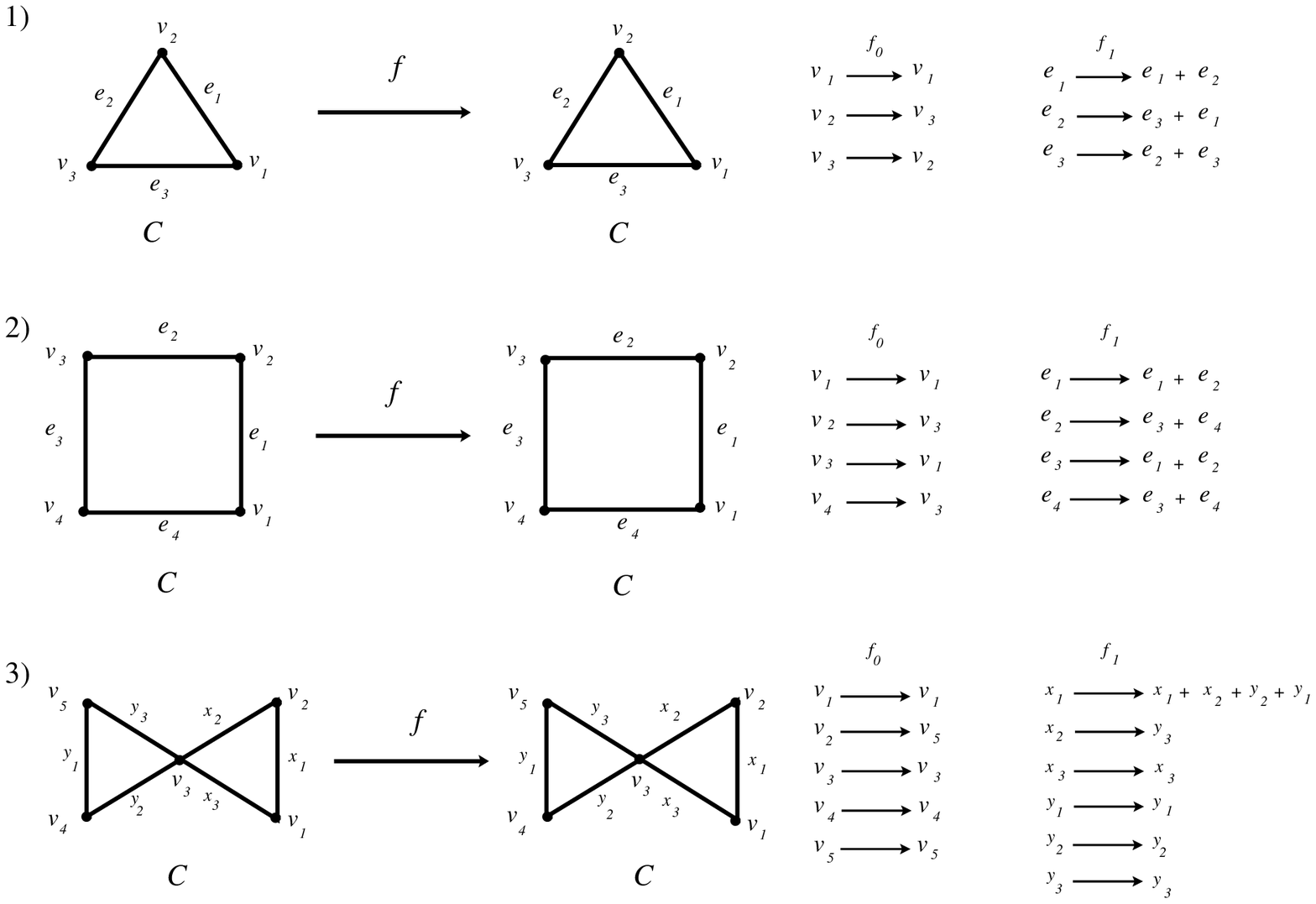}}
\medskip
\centerline{\small Figure 2. Examples of quasi-isomorphisms.}
\bigskip

The chain complexes in the examples are chain complexes over the
real field $\mathbb R$.

Example 1. We compute the torsion of $f:C \rightarrow C$ as {\rm
1)} in Figure 2 by definition. Suppose that $C_0=\langle v_1, v_2,
v_3 \rangle$, $C_1=\langle e_1, e_2, e_3 \rangle$, $c_0=(v_1, v_2,
v_3)$, and $c_1=(e_1, e_2, e_3)$ and the $0$-th boundary map
$\partial_0:C_1 \rightarrow C_0$ is defined by
$\partial_0(e_1)=v_2-v_1$, $\partial_0(e_2)=v_3-v_2$, and
$\partial_0(e_3)=v_1-v_3$. Then we can easily show that $f$ is a
chain map. Since $B_0(C)={\rm sp}(v_2-v_1, v_3-v_2, v_1-v_3)={\rm
sp}(v_2-v_1, v_3-v_2)$, we can take $b_0=(v_2-v_1, v_3-v_2)$.
Also, $H_0(C)=C_0/B_0(C)$. To take a basis $h_0$ for $H_0(C)$, we
need a vector in $C_0$ not contained in $B_0(C)$. For this, we
take an orthogonal vector $v_1+v_2+v_3$ to $B_0(C)$, so we take
$h_0=(v_1+v_2+v_3+B_0(C))$. Since
$f_{0*}(v_1+v_2+v_3+B_0(C))=v_1+v_3+v_2+B_0(C)$, $f_{0*}$ is the
identity map. Suppose that $\partial_0(r_1e_1 + r_2e_2 +
r_3e_3)=0$. Then $r_1(v_2-v_1) + r_2(v_3-v_2) + r_3(v_1-v_3)=0$.
Hence, $r_1=r_2=r_3$ and we have $Z_0(C)=\langle e_1+e_2+e_3
\rangle$. Since $B_0(C)=0$, $H_0(C)=Z_0(C)$. Take
$h_1=(e_1+e_2+e_3+0)$. Since
$f_{1*}(e_1+e_2+e_3+0)=2(e_1+e_2+e_3)+0$ and $2(e_1+e_2+e_3)
\notin B_1(C)=0$, $f_{1*}$ is an isomorphism.

Now, we compute the torsion $\tau(f)$ of $f$. Notice that
$b_{-1}=b_1=\emptyset$. Since $(b_0h_0)b_{-1}=(v_2-v_1, v_3-v_2,
v_1+v_2+v_3)$, we have $[(b_0h_0)b_{-1}/c_0]={\rm
det}\begin{pmatrix} -1 & 1 & 0
\\ 0 & -1 & 1 \\ 1 & 1 & 1 \end{pmatrix}=3$. Also, $(b_0f_{0*}(h_0))b_{-1}=
(v_2-v_1, v_3-v_2, v_1+v_3+v_2)$, hence,
$[(b_0f_{0*}(h_0))b_{-1}/c_0]=3$. Similarly, since
$(b_1h_1)b_0=(e_1+e_2+e_3, e_1, e_2)$ and
$(b_1f_{1*}(h_1))b_0=(2(e_1+e_2+e_3), e_1, e_2)$, we have
$[(b_1h_1)b_0/c_1]=1$ and $[(b_1f_{1*}(h_1))b_0/c_1]=2$.
Therefore,
$$\tau(f)=\left(\frac{[(b_0h_0)b_{-1}/c_0]}{[(b_0f_{0*}(h_0))b_{-1}/c_0]}
\right)^{-1}\left(\frac{[(b_1h_1)b_0/c_1]}{[(b_1f_{1*}(h_1))b_0/c_1]}
\right)=\left(\frac{3}{3}\right)^{-1}\left(\frac{1}{2}\right)=\frac{1}{2}\,.$$

Example 2. We use Theorem 3.20 to compute the torsion of $f:C
\rightarrow C$ as {\rm 2)} in Figure 2. Suppose that $C_0=\langle
v_1, v_2, v_3, v_4 \rangle$, $C_1=\langle e_1, e_2, e_3, e_4
\rangle$, $c_0=(v_1, v_2, v_3, v_4)$, and $c_1=(e_1, e_2, e_3,
e_4)$ and the $0$-th boundary map $\partial_0:C_1 \rightarrow C_0$
is defined by $\partial_0(e_1)=v_2-v_1$,
$\partial_0(e_2)=v_3-v_2$, $\partial_0(e_3)=v_4-v_3$, and
$\partial_0(e_4)=v_1-v_4$. Then $f$ is a chain map. Since
$B_0(C)={\rm sp}(v_2-v_1, v_3-v_2, v_4-v_3, v_1-v_4)={\rm
sp}(v_2-v_1, v_3-v_2, v_4-v_3)$, we can take $b_0=(v_2-v_1,
v_3-v_2, v_4-v_3)$. Also, $H_0(C)=C_0/B_0(C)$. To take a basis
$h_0$ for $H_0(C)$, we need a vector in $C_0$ not contained in
$B_0(C)$. As in Example 1, we can take
$h_0=(v_1+v_2+v_3+v_4+B_0(C))$. Since
$f_{0*}(v_1+v_2+v_3+v_4+B_0(C))=v_1+v_3+v_1+v_3+B_0(C)$ and
$(v_1+v_2+v_3+v_4)-(v_1+v_3+v_1+v_3)=(v_2-v_1)+(v_4-v_3) \in
B_0(C)$, $f_{0*}$ is the identity map. Hence, ${\rm
det}\,f_{0*}=1$. As in Example 1, we have $Z_0(C)=\langle
e_1+e_2+e_3+e_4 \rangle$. Let us take $h_1=(e_1+e_2+e_3+e_4+0)$.
Then $f_{1*}(e_1+e_2+e_3+e_4+0)=2(e_1+e_2+e_3+e_4)+0$. Since
$B_1(C)=0$, $2(e_1+e_2+e_3+e_4) \notin B_1(C)$. Hence, $f_{1*}$ is
an isomorphism and ${\rm det}\,f_{1*}=2$. Therefore, by Theorem
3.20, $$\tau(f)=\left(\frac{1}{{\rm
det}\,f_{0*}}\right)^{-1}\left(\frac{1}{{\rm det}\,f_{1*}}\right)
=\left(\frac{1}{1}\right)^{-1}\left(\frac{1}{2}\right)=\frac{1}{2}\,.$$

First two examples are about 2-fold covering maps. Finally, we try
to get the torsion of a little bit more complicated
quasi-isomorphism.

Example 3. We also use Theorem 3.20 for the torsion of $f:C
\rightarrow C$ as {\rm 3)} in Figure 2. Suppose that $C_0=\langle
v_1, v_2, v_3, v_4, v_5 \rangle$, $C_1=\langle x_1, x_2, x_3, y_1,
y_2, y_3 \rangle$, $c_0=(v_1, v_2, v_3, v_4, v_5)$, and $c_1=(x_1,
x_2, x_3, y_1, y_2, y_3)$ and the $0$-th boundary map
$\partial_0:C_1 \rightarrow C_0$ is defined by
$\partial_0(x_1)=v_2-v_1$, $\partial_0(x_2)=v_3-v_2$,
$\partial_0(x_3)=v_1-v_3$, $\partial_0(y_1)=v_5-v_4$,
$\partial_0(y_2)=v_4-v_3$, and $\partial_0(y_3)=v_3-v_5$. Then $f$
is a chain map. Since $B_0(C)={\rm sp}(v_2-v_1, v_3-v_2, v_1-v_3,
v_5-v_4, v_4-v_3, v_3-v_5)$, we can take $b_0=(v_2-v_1, v_3-v_2,
v_5-v_4, v_4-v_3)$. Also, $H_0(C)=C_0/B_0(C)$. As in Example 1, we
take $h_0=(v_1+v_2+v_3+v_4+v_5+B_0(C))$. Since
$f_{0*}(v_1+v_2+v_3+v_4+v_5+B_0(C))=v_1+v_3+v_4+2v_5+B_0(C)$ and
$(v_1+v_3+v_4+2v_5)-(v_1+v_2+v_3+v_4+v_5)=v_5-v_2=(v_5-v_4)+(v_4-v_3)
+(v_3-v_2) \in B_0(C)$, $f_{0*}$ is the identity map. Hence, ${\rm
det}\,f_{0*}=1$.

Suppose that $\partial_0(r_1x_1 + r_2x_2 + r_3x_3 + s_1y_1 +
s_2y_2 + s_3y_3)=0$. Then $r_1(v_2-v_1) + r_2(v_3-v_2) +
r_3(v_1-v_3)+s_1(v_5-v_4) + s_2(v_4-v_3) + s_3(v_3-v_5)=0$ and we
have $r_1=r_2=r_3$ and $s_1=s_2=s_3$. Hence, $Z_0(C)=\langle
x_1+x_2+x_3, y_1+y_2+y_3 \rangle$. Since $B_0(C)=0$,
$H_0(C)=Z_0(C)$. Take $h_1=(x_1+x_2+x_3+0, y_1+y_2+y_3+0)$. Since
$f_{1*}(x_1+x_2+x_3+0)=(x_1+x_2+x_3+0)+(y_1+y_2+y_3+0)$ and
$f_{1*}(y_1+y_2+y_3+0)=y_1+y_2+y_3+0$. Hence, the matrix
representation of $f_{1*}$ for $h_1$ is $\begin{pmatrix} 1 & 1 \\
0 & 1 \end{pmatrix}$ and ${\rm det}\,f_{1*}=1 \ne 0$. Therefore
$f_{1*}$ is an isomorphism and, by Theorem 3.20, we have
$$\tau(f)=\left(\frac{1}{{\rm
det}\,f_{0*}}\right)^{-1}\left(\frac{1}{{\rm det}\,f_{1*}}\right)
=\left(\frac{1}{1}\right)^{-1}\left(\frac{1}{1}\right)=1.$$
\\

\end{document}